\documentclass[11pt]{amsart}

\usepackage{amsfonts, amsmath, amssymb,amsthm}

\newcommand{\vanish}[1]{}

\newcommand{\Z}{\mathbb{Z}}  
\newcommand{\N}{\mathbb{N}}  
\newcommand{\Q}{\mathbb{Q}}  
\newcommand{\R}{\mathbb{R}}  
\newcommand{\PP}{\mathbb{P}}  

\newcommand{\QQ}{\mathcal{Q}}

\newcommand{\qs}{QSym}
\newcommand{\I}{\mathcal{I}}
\newcommand{\T}{\mathcal{T}}
\newcommand{\cd}{{\bf cd}}
\newcommand{\bc}{{\bf c}}
\newcommand{\bd}{{\bf d}}
\newcommand{\ba}{{\bf a}}
\newcommand{\bb}{{\bf b}}

\providecommand{\abs}[1]{\lvert#1\rvert}

\newcommand{\rk}{\rho}  
  
\newcommand{\cal}{\mathcal} 
\newcommand{\K}{{\mathcal K}}

\newtheorem{theorem}{Theorem}
\newtheorem{proposition}[theorem]{Proposition}

\newtheorem{corollary}[theorem]{Corollary}

\theoremstyle{definition}

\newtheorem{example}[theorem]{Example}
\newtheorem{conjecture}[theorem]{Conjecture}
                                                                              
\theoremstyle{remark}
\newtheorem{remark}[theorem]{Remark} 

\numberwithin{theorem}{section}
\numberwithin{equation}{section}
                                                                               
\begin{document}

\title[Kazhdan-Lusztig polynomials]{Quasisymmetric functions and
Kazhdan-Lusztig polynomials}

\author{Louis J.\ Billera}
\address{Department of Mathematics, Cornell University, Ithaca, NY
14853-4201}
\email{billera@math.cornell.edu}
\thanks{The first author was supported in part by NSF grants
DMS-0100323 and DMS-0555268.
The second author was partially supported by EU grant CHRT-CT-2001-00272.
This work was begun while both authors enjoyed the hospitality of
the Mittag-Leffler Institute, Djursholm, Sweden.}

\author{Francesco Brenti}
\address{Dipartimento di Matematica, Universit\'a di Roma ``Tor Vergata'',
Via della Ricerca Scientifica, 1, I-00133 Roma, Italy}
\email{brenti@mat.uniroma2.it}


\begin{abstract}
We associate a quasisymmetric function to any Bruhat interval in a  general Coxeter group.  This
association can be seen to be a morphism of Hopf algebras to the subalgebra of all peak functions,
leading to an extension of the {\bf cd}-index of convex polytopes.
We show how the Kazhdan-Lusztig polynomial of the Bruhat interval can be expressed in
terms of this complete {\bf cd}-index and otherwise explicit combinatorially defined polynomials.
In particular, we obtain the simplest closed 
formula for the Kazhdan-Lusztig polynomials that holds in complete 
generality.
\end{abstract}

\maketitle

\setcounter{tocdepth}{1}
\tableofcontents 


\section{Introduction}
The Kazhdan-Lusztig polynomials of a Coxeter group are of fundamental importance in
representation theory and in the geometry and topology of Schubert varieties. Defined by
means of two separate recursions, they have proved difficult to unravel in any
straightforward manner. Here, we reduce their computation to the computation of
another invariant of Coxeter groups.  This new
invariant, a quasisymmetric function that can be encoded into a noncommutative polynomial
in two variables that we call the \emph{complete {\bf cd}-index}, has interesting algebraic
and combinatorial properties.
We express the Kazhdan-Lusztig polynomial of any Bruhat interval in terms
of its complete {\bf cd}-index and otherwise explicit combinatorially defined polynomials.
In particular, we obtain the simplest closed 
formula for the Kazhdan-Lusztig polynomials that holds in complete 
generality.

The complete {\bf cd}-index is defined by means of a quasisysymmetric function
$\widetilde{F}(u,v)$ associated to every Bruhat interval $[u,v]$ in a Coxeter group $W$.
The association $[u,v] \mapsto \widetilde{F}(u,v)$  can be viewed as a morphism of
Hopf algebras, suggesting the possibility of a filtered version of the theory of combinatorial
Hopf algebras \cite{ABS}, which predicts the existence of graded maps to the quasisymmetric
functions in general combinatorial settings.

The theory of Kazhdan-Lusztig polynomials of Weyl groups is analogous to that
of the $g$-polynomials of rational convex polytopes in that they both compute the local
intersection cohomology of certain varieties (Schubert and toric, respectively)
associated to these objects \cite{KL2}, \cite{StaH}. Also, the recursions that define them
in a general Coxeter group (respectively, Eulerian partially ordered set) have
the same form. The $g$-polynomial of an Eulerian partially ordered set is known
to depend only on the number of chains of certain types. These flag numbers are
most succintly represented by the $\cd$-index. On the other hand, for the
Kazhdan-Lusztig polynomials it is not even clear that there can be a completely combinatorial
description (see, {\it e.g.}, \cite[\S 5.6]{BjBr}), so there is no straightforward way to
generalize this
combinatorial description of the $g$-polynomials to the Kazhdan-Lusztig polynomials.
Yet, the quasisymmetric functions $\widetilde{F}(u,v)$ do seem to capture much
of the spirit of the combinatorial setting of the $g$-polynomials, without themselves
being obviously combinatorial.

In the remainder of this section we give the necessary background in Coxeter groups and
Kazhdan-Lusztig polynomials and in the use of quasisymmetric functions in poset
enumeration.  In Section \ref{sect2} we introduce the $R$-quasisymmetric function of a Bruhat
interval, define from this the complete {\bf cd}-index of the interval and give some  of its algebraic
properties.
Section \ref{sect3} gives an expression for the Kazhdan-Lusztig polynomials
in terms of the complete {\bf cd}-index and certain lattice path enumerators.
Using this expression, we give
in Section \ref{sect4} the representation of
the Kazhdan-Lusztig polynomials in terms of an explicit polynomial basis defined
by means of  the ballot polynomials.  The coefficients of this representation are
given as linear forms in the complete $\cd$-index.  We show by an example
that no such representation exists in terms of the ordinary $\cd$-index alone; that is,
the Kazhdan-Lusztig polynomials of a Bruhat interval can not be calculated, in general,
from the flag $f$-vector of the underlying Eulerian poset.
In Section \ref{sect5} we give a formula for computing the complete
{\bf cd}-index of any Bruhat interval in terms of explicit combinatorial
quantities associated to the Bruhat interval.
Finally, in Section \ref{sect7} we conjecture nonnegativity of the complete \cd-index
and show that one consequence of this conjecture holds in the case of finite
Coxeter and affine Weyl groups.

\subsection{Coxeter groups and Kazhdan-Lusztig polynomials}
We follow \cite{Hum} for general Coxeter groups notation
and terminology.
In particular, given
a Coxeter system $(W,S)$ and $u \in W$ we denote by $l
(u )$ the length of $u $ in $W$, with respect to $S$.
We denote by $e$
the identity of $W$, and we let $T \overset{\rm def}{=}
\{ u s u ^{-1} : u \in W, \; s \in S \}$ be the
set of \emph{reflections} of $W$.
We will
always assume that $W$ is partially ordered by {\em Bruhat
order}. Recall (see, {\it e.g.}, \cite[\S 5.9]{Hum}) that this means that
$x \leq y$ if and only if there exist $r \in \N$ and $t_{1},
\ldots , t_{r} \in T$ such that $ t_{r} \cdots  t_{1} \, x=y$ and
$l (t_{i} \cdots t_{1} \, x) > l(t_{i-1} \cdots t_{1}x)$ for $i=1,
\ldots ,r$.
Given $u,v \in W$ we let $[u,v] \overset{\rm def}{=} \{ x \in W
: u \leq x \leq v \}$.  We consider $[u,v]$ as a poset with the 
partial ordering induced by $W$.
It is well known (see, {\it e.g.}, 
\cite{BjBr}, Corollary 2.7.11)
that intervals of $W$ (and their duals) are Eulerian
posets.

Let $A \subseteq T$ and $W'$ be the subgroup of $W$ generated by
$A$. Following \cite[\S 8.2]{Hum},
we call $W'$ a {\em reflection subgroup} of $W$. It is then
known (see, {\it e.g.}, \cite{Hum}, Theorem 8.2) that
$(W',S')$ is again a Coxeter system where
$S' \overset{\rm def}{=} \{ t \in T: \; N(t) \cap W' = \{ t \} \} $,
and $ N(w) \overset{\rm def}{=} \{ t \in T: \; l(wt)<l(w) \}$.
We say that
$W'$ is a {\em dihedral reflection subgroup} if $|S'|=2$ ({\it i.e.},
if $(W',S')$ is a dihedral Coxeter system). Following \cite{Dy}
we say that a total ordering $<_{T}$ of $T$ is a {\em
reflection ordering} if, for any dihedral reflection
subgroup $W'$ of $W$, we have that either $a <_{T} aba <_{T}
ababa <_{T} \cdots <_{T} babab <_{T}
bab <_{T}  b$ or $b <_{T} bab <_{T} babab <_{T} \cdots
<_{T} ababa <_{T} aba <_{T} a$ where
$\{ a,b \} \overset{\rm def}{=} S'$.
 The existence of reflection orderings
(and many of their properties) is proved in \cite[\S 2]{Dy}
(see also \cite[\S 5.2]{BjBr}).
Throughout this work we will always assume that we have fixed (once and
for all) a reflection ordering $<_{T}$ of $T$.

We denote by ${\mathcal H}(W)$ the {\em Hecke algebra} associated
to $W$. Recall (see, {\it e.g.}, \cite[Chap. 7]{Hum})
  that this is the free $\Z[q,q^{-1}]$-module
having the set $\{ T_{w}: \; w \in W \}$ as a basis and
multiplication such that
\begin{equation}
\label{2.1.2}
T_{w} T_{s} = \left\{ \begin{array}{ll}
T_{ws}, & \mbox{if $l(ws) > l(w)$,} \\
q T_{ws} +(q-1)T_{w}, & \mbox{if $l(ws) < l(w)$,}
\end{array}
\right. 
\end{equation}
for all $w \in W$ and $s \in S$. It is well known that this is an
associative algebra having $T_{e}$ as unity and that each basis
element is invertible in ${\cal H} (W)$. More precisely, we have
the following result (see \cite[Proposition 7.4]{Hum}).
\begin{proposition}
\label{2.2}
Let $v \in W$. Then
\[ (T_{v^{-1}})^{-1} =q^{-l(v)} \, \sum _{u \leq v} (-1)^{l(v)-l(u)}
\, R_{u,v}(q) \, T_{u} \, , \]
where $R_{u,v}(q) \in \Z[q]$.
\end{proposition}
The polynomials $R_{u,v}(q)$ defined by the previous proposition
are called the {\em $R$-polynomials } of $W$.
It is easy to see that $R_{u,v}(q)$ is a monic polynomial of degree 
$l(u,v) \overset{\rm def}{=} l(v)-l(u)$,
and that $R_{u,u}(q)=1$, for all $u,v \in W$, $u \leq v$.
It is customary to let $R_{u,v}(q) \overset{\rm def}{=}
0$ if $u \not \leq v$.
We then  have the 
following fundamental result that follows from (\ref{2.1.2})
and Proposition \ref{2.2}
 (see \cite[\S 7.5]{Hum}).
\begin{theorem}
\label{Rpoly}
Let $u,v \in W$ and $s \in S$ be such that $l(vs)<l(v)$. Then
\begin{equation*}
 R_{u,v}(q)= \left\{ \begin{array}{ll}
R_{us,vs}(q), & \mbox{if $l(us)<l(u)$,} \\
q R_{us, vs}(q)+(q-1)R_{u,vs}(q), & \mbox{if $l(us)>l(u)$.}
\end{array} \right. 
\end{equation*}
\end{theorem}
Note that the preceding theorem can be used to inductively 
compute the $R$-polynomials. 
Theorem \ref{Rpoly} has also the following simple  but important consequence
(see, for example, \cite[Proposition 5.3.1]{BjBr}).
\begin{proposition}
Let $u, v \in W$. Then there exists a (necessarily 
unique) polynomial $\widetilde{R}_{u,v}(q) \in \N [q]$ such that
\begin{equation*}
 R_{u,v}(q) =q^{\frac{1}{2}(l(v)-l(u))} \, \widetilde{R}_{u,v} \left( q^{
\frac{1}{2}} - q^{- \frac{1}{2}} \right) \, . 
\end{equation*}
\label{Rtilde}
\end{proposition}
Combinatorial interpretations of the coefficients of
$\widetilde{R}_{u,v}(q)$ have been given by V. Deodhar \cite{Deo}
and by M. Dyer \cite{Dy} (see \cite[Theorem 5.3.7]{BjBr} and \cite[Theorem 5.3.4]{BjBr}).

The $R$-polynomials can be
used to define
the Kazhdan-Lusztig polynomials.
The following result is not hard to prove (and, in fact, holds in
much greater generality, see \cite[Corollary 6.7]{StaLH} and
\cite[Example 6.9]{StaLH}) and a proof can be found, {\it e.g.}, in \cite[\S 7.9--11]{Hum}
or \cite[\S 2.2]{K-L}.
\begin{theorem}
\label{KLdefn}
There is a unique family of polynomials $\{ P_{u,v}(q) \}
_{u,v \in  W}  \subseteq \Z [q]$, such that, for all
$u,v \in W$,
\begin{enumerate}
\item[1.]
$P_{u,v}(q)=0$ if $u \not \leq v$;
\item[2.]
$P_{u,u}(q)=1$;
\item[3.]
deg$(P_{u,v}(q)) < \frac{1}{2}\left(
l(v)-l(u) \right) $, if $u < v$;
\item[4.]
\begin{equation*}
q^{l(v)-l(u)} \, P_{u,v} \left( \frac{1}{q} \right)
= \sum _{u \leq z \leq v}   R_{u,z}(q) \,
P_{z,v}(q) \, ,
\end{equation*}
if $u \leq v$.
\end{enumerate}
\end{theorem}
The polynomials $P_{u,v}(q)$ defined by the preceding
theorem are called the
{\em Kazhdan-Lusztig polynomials} of $W$.

\subsection{Paths in Bruhat graphs.}
Recall that a {\em composition} of a positive integer $n$ is
a finite sequence of positive integers $\alpha = \alpha_1 \cdots  \alpha_k$
such that $\sum_i \alpha_i = n$. In this case we write $\alpha \models n$
and we define $l(\alpha) = k$, $|\alpha| = n$, and
$\alpha ^{\ast} = \alpha_k \cdots \alpha_1$.
(There should be no confusion with our using the same notation $l(\cdot)$
for the lengths of
compositions and elements of $W$.)
Given two compositions $\alpha_{1} \cdots \alpha_{s},\beta_{1} \cdots \beta_{t}$
of $n$ we say that $\alpha_{1} \cdots \alpha_{s}$ {\em refines}
$\beta_{1} \cdots \beta_{t}$ if there exist $1\leq i_{1} < i_{2}
< \cdots < i_{t-1} < s$ such that $\sum_{j=i_{k-1}+1}^{i_{k}}
\alpha_{j} = \beta_{k}$ for $k=1,\ldots,t$ (where $i_{0}
\overset{\rm def}{=} 0$, $i_{t} \overset{\rm def}{=}
s$). We then write $\alpha_{1} \cdots \alpha_{s} \preceq
\beta_{1} \cdots \beta_{t}$. We denote by $C$ the set of all finite sequences
of positive integers ({\it i.e.}, the set of all compositions).

Recall (see \cite{Hum}, \S 8.6, or \cite{Dy1})
that the {\em Bruhat graph} of a Coxeter system $(W,S)$
is the directed graph $B(W,S)$ obtained by taking $W$ as vertex
set and putting a directed edge from $x$ to $y$ if and only if
$y x^{-1} \in T$ and $l(x)<l(y)$.
We call the directed paths of $B(W,S)$ {\em Bruhat paths}.
These should be distinguished from \emph{chains} in the Bruhat
\emph{order}.  The vertex set of a Bruhat path is always a chain, but
not all chains form Bruhat paths.

Given a Bruhat path $\Delta = (a_{0}, a_{1} , \ldots ,
a_{r})$ from $a_{0}$ to $a_{r}$, we define
 its {\em length} to be $l(\Delta) \overset{\rm def}{=} r$,
its
 {\em descent set}, with respect to $<_{T}$, to be
\[ D(\Delta ) \overset{\rm def}{=} \{ i \in [r-1]: \; a_{i}(a_{i-1})^{-1}
>_{T} a_{i+1} (a_{i})^{-1} \} \, , \]
and  its {\em descent composition}
to be the unique composition
${\cal D}(\Delta )$ of $r$ corresponding to $D(\Delta )$ under
the usual bijection between compositions and sets
$\beta_{1}\beta_{2}\cdots \beta_{k} \mapsto
\{\beta_{1},\beta_{1}+\beta_{2},\dots, \beta_{1}+\cdots + \beta_{k-1}\}$.
We will denote the inverse of this bijection by $co(\cdot)$.
Given $u,v \in W$, and $k \in \N$,
 we denote by $B_{k}(u,v)$ the set of all directed
paths in $B(W,S)$ from $u$ to $v$ of length $k$, and we let $B(u,v)
\overset{\rm def}{=} \bigcup _{k \geq 0} B_{k}(u,v)$.
For $u,v \in W$, and $\alpha \in C$,
we let, following \cite{BreLMS},
\begin{equation}
\label{}
c _{\alpha } (u,v) \overset{\rm def}{=}
| \{ \Delta  \in B_{|\alpha |}(u,v) : \;
{\cal D}(\Delta ) \succeq \alpha \} | , 
\end{equation}
and
\begin{equation}
\label{balpha}
 b_{\alpha }(u,v) \overset{\rm def}{=}
| \{ \Delta  \in B_{|\alpha |}(u,v): \; {\cal D}(\Delta )
= \alpha \} | .
\end{equation}
 Note that these definitions  imply
that
\begin{equation}
\label{2.IE.CB}
 c_{\alpha }(u,v) = \sum _{\{ \beta \models n : \beta \succeq \alpha \} }
b_{\beta}(u,v) 
\end{equation}
for all $u,v \in W$ and $\alpha \models n$ ($n \in \PP$).

The following result follows from \cite[Proposition 4.4]{BreLMS}.
 Given a polynomial $P(q)$, and $i \in \Z$, we
denote by $[q^{i}](P)$ the coefficient of $q^{i}$ in $P(q)$.
\begin{proposition}
\label{2.LMS}
Let $u,v \in W$, $u \leq v$, and $\alpha \in C$. Then
\[ c_{\alpha } (u,v) = \sum _{(u_{0}, \ldots , u_{r}) \in
 C_{r}(u,v)}\ \prod _{j=1}^{r} [q^{\alpha _{j}}] (\widetilde{R}_{u_{j-1}, u_{j}}) \]
where $C_{r}(u,v)$ denotes the set of all chains of length $r$
(totally ordered subsets of cardinality $r+1$)
 from $u$ to $v$, and $r \overset{\rm def}{=} l(\alpha )$.
\end{proposition}
Proposition \ref{2.LMS} shows, in particular, that $c_{\alpha}(u,v)$
(and hence $b_{\alpha }(u,v)$) are
independent of the total reflection ordering $<_{T}$ used to define them.

For $j \in \Q$ 
we define an operator $D_{j} : 
\R [q] \rightarrow \R [q]$  by letting 
\[ D_{j} \left( \sum _{i \geq 0} a_{i}q^{i} \right) \overset{\rm def}{=}
\sum _{i =0}^{\lfloor j \rfloor} a_{i}q^{i} . \]
For $\alpha \in C$ we define, following \cite{BrJAMS}, a polynomial
$\Psi_{\alpha }(q) \in \Z[q]$ inductively as follows,
\begin{equation*}
\Psi_{\alpha }(q) \overset{\rm def}{=} (q-1)^{\alpha _{1}}
D_{\frac{|{\alpha^{-} }|-1}{2}} (\Psi_{{\alpha^{-}}}(q))
\end{equation*}
if $l(\alpha ) \geq 2$, and
\begin{equation*}
\Psi_{\alpha }(q) \overset{\rm def}{=} (q-1)^{|\alpha |} ,
\end{equation*} 
 if $l(\alpha )=1$.  Here $\alpha^{-}= \alpha_{2} \cdots \alpha_{k}$ when
 $\alpha= \alpha_{1} \cdots \alpha_{k}$.
 For $ n \in  \PP$ and  $\beta \models n$ we then let
 \begin{equation*}
  \Upsilon _{\beta }(q) \overset{\rm def}{=}
  \sum _{\{ \alpha \models n : \; \alpha \preceq \beta \} } (-1)^{l(\alpha )}
 \Psi_{\alpha } (q) .
 \end{equation*}
We then have the following result, whose proof can be found in \cite[Theorem 5.5.7]{BjBr}. 
\begin{theorem}
\label{1.8}
Let $u,v \in W$, $u<v$.  Then
 \begin{equation*}
 P_{u,v}(q) - q^{l(u,v)} P_{u,v} \left( \frac{1}{q} \right) =
\sum _{\beta \in C}  q^{\frac{l(u,v)- |\beta | }{2} }
\Upsilon _{\beta }(q) \, b_{\beta } (u,v) .
\end{equation*}
\end{theorem}

Let $(W_{1},S_{1})$ and $(W_{2},S_{2})$ be two Coxeter systems.  It is clear that the
direct product $W = W_{1}\times W_{2}$ is again a Coxeter group with respect to the
generating set $S=S_{1} \sqcup S_{2}$ (disjoint union), having Dynkin diagram the disjoint
union of the diagrams of $(W_{1},S_{1})$ and $(W_{2},S_{2})$.  It then follows immediately from the subword property of Bruhat order (see, {\it e.g.}, \cite[Theorem 2.2.2]{BjBr}) that $W$, as a poset
under Bruhat order, is isomorphic to the direct product of $W_{1}$ and $W_{2}$ as posets under
Bruhat order.  Thus $(u_{1}, u_{2}) \le (v_{1}, v_{2})$ in $W$ if and only if $u_{1}\le v_{1}$ in
$W_{1}$ and $u_{2}\le v_{2}$ in $W_{2}$.  Thus Theorems \ref{Rpoly} and \ref{KLdefn} imply
the following multiplicative formulas.  We include a proof for lack of an adequate reference.

\begin{proposition}
\label{productform}
Let $u_{1}, v_{1} \in W_{1}$ and  $u_{2}, v_{2} \in W_{2}$.  Then
\begin{align}
\label{R.R}
R_{(u_{1},u_{2}),(v_{1},v_{2})} &= R_{u_{1}, v_{1}} \cdot R_{u_{2},v_{2} }, \text{~and}\\
\label{P.P}
P_{(u_{1},u_{2}),(v_{1},v_{2})} &= P_{u_{1}, v_{1}} \cdot P_{u_{2},v_{2} }.
\end{align}
\begin{proof}
First note that by the comments preceding this proposition, \eqref{R.R} and \eqref{P.P} hold
if $(u_{1}, u_{2}) \not\le (v_{1}, v_{2})$, so assume that $(u_{1}, u_{2}) \le (v_{1}, v_{2})$.
We first prove \eqref{R.R} by induction on $l(v_{2})$.If $v_{2}=e$ then $u_{2}=e$ and \eqref{R.R}
follows from Theorem \ref{Rpoly} and the fact that $D((u,e))=
\{ (s,e) \, : s \in D(u) \}$ for all $u\in W_{1}$.  (Note that
$l_{W}\big((u_{1},u_{2}),(v_{1},v_{2})\big)=l_{W_{1}}(u_{1},v_{1})+l_{W_{2}}(u_{2},v_{2})$.)
If $v_{2}>e$ and $s \in D(v_{2})$, then $(e,s)\in D((v_{1},v_{2}))$ and therefore, by Theorem
\ref{Rpoly} and our induction hypothesis,
\begin{align*}
R_{(u_{1},u_{2}),(v_{1},v_{2})} &= q\, R_{(u_{1}, u_{2}s),(v_{1}, v_{2}s)} 
					+ (q-1)\, R_{(u_{1},u_{2}),(v_{1},v_{2}s)}\\
		&= R_{u_{1},v_{1}}\, \left(q\, R_{u_{2}s, v_{2}s} + (q-1) \, R_{u_{2}, v_{2}s}\right)\\
		&= R_{u_{1},v_{1}} \cdot R_{u_{2},v_{2}}
\end{align*}
if $s \notin D(u_{2})$ (so $(e,s) \notin D((u_{1},u_{2})))$, while
\begin{align*}
R_{(u_{1},u_{2}),(v_{1},v_{2})} = R_{(u_{1},u_{2}s),(v_{1},v_{2}s)} 
			&= R_{u_{1},v_{1}} \; R_{u_{2}s,v_{2}s}\\
			&= R_{u_{1},v_{1}} R_{u_{2},v_{2}}
\end{align*}
if $s\in D(u_{2})$, as desired.

We conclude by proving \eqref{P.P} by induction on $l_{W}\big( (u_{1},u_{2}), (v_{1},v_{2}) \big)$.
The result is clear if $l_{W}\big( (u_{1},u_{2}), (v_{1},v_{2}) \big)=0$.  If
$l_{W}\big( (u_{1},u_{2}), (v_{1},v_{2}) \big)>0$ then by Theorem \ref{KLdefn}, our induction
hypothesis, \eqref{R.R} and the comments preceding this proposition, we have
\begin{align*}
q&^{l( (u_{1},u_{2}), (v_{1},v_{2}))}\, P_{(u_{1},u_{2}),(v_{1},v_{2})}\left(\frac 1 q \right)  -
				P_{(u_{1},u_{2}),(v_{1},v_{2})}(q)  \\
				&=\sum_{(u_{1},u_{2}) < (x_{1},x_{2}) \le (v_{1},v_{2})}
					R_{(u_{1},u_{2}),(x_{1},x_{2})}(q) \; P_{(x_{1},x_{2}),(v_{1},v_{2})}(q) \\
				&= \sum_{u_{1}\le x_{1}\le v_{1}} R_{u_{1},x_{1}}(q)\, P_{x_{1},v_{1}}(q)
				\, \sum_{u_{2}\le x_{2}\le v_{2}} R_{u_{2},x_{2}}(q)\, P_{x_{2},v_{2}}(q)\\
				& \hskip13.5pc - P_{u_{1},v_{1}}(q) \; P_{u_{2},v_{2}}(q)\\
				&= q^{l(u_{1},v_{1})}\; P_{u_{1},v_{1}}\left(\frac 1 q \right)\, 
					 q^{l(u_{2},v_{2})}\; P_{u_{2},v_{2}}\left(\frac 1 q\right)
					 - P_{u_{1},v_{1}}(q) \; P_{u_{2},v_{2}}(q),
\end{align*}
and \eqref{P.P} follows.
\end{proof}

\end{proposition}

Throughout this work (unless otherwise explicitly stated) 
$(W,S)$ denotes a fixed (but arbitrary) Coxeter system.

\subsection{Quasisymmetric functions and poset enumeration}
A \emph{quasisymmetric function} is a formal power series in countably many
variables that has bounded degree and whose coefficients are invariant
under shifts of the variables that respect their order.
We assume here that the reader is familiar with the basics of the theory of
quasisymmetric functions, for example, as described in \cite[\S 7.19]{ECII}.
We denote by $\qs\subset \Q[[x_{1},x_{2}, \dots ]]$ the algebra of all quasisymmetric
functions (with rational coefficients).  $\qs$ is a graded algebra with the
usual grading of power series; we denote by $\QQ_{i}$ the
i$^{th}$ homogeneous part of $\qs$ and so
$$\qs = \QQ_{0}\oplus \QQ_{1} \oplus \cdots .$$

In particular, we will make use of the \emph{monomial basis} $\{ M_{\alpha}\}_{\alpha \in C}$ and
the \emph{fundamental basis} $\{L_{\alpha}\}_{\alpha \in C}$ for $\qs$, where for a composition
$\alpha = \alpha_{1}\alpha_{2}\cdots\alpha_{k}$, $k>0$, $\alpha_{i}>0$,
$$M_{\alpha} = \sum_{i_{1}<i_{2}<\cdots <i_{k}}  x_{i_{1}}^{\alpha_{1}} x_{i_{2}}^{\alpha_{2}} \cdots
x_{i_{k}}^{\alpha_{k}} $$
and
\begin{equation}
L_{\beta}=  \sum _{\{ \alpha \models |\beta| : \; \alpha \preceq \beta \} } M_{\alpha}.
\label{MtoL}
\end{equation}
We include as well the empty composition $\alpha = {\bf 0}$ (the case $k=0$); here we set
$M_{\bf 0}= L_{\bf 0} = 1$. 
Note that the degree of $M_{\alpha}$ and $L_{\alpha}$ is $|\alpha|$.
Occasionally, it will be useful to index
these bases by subsets of $[n-1]$ instead of compositions of $n$, using the standard bijection
between compositions and subsets already mentioned.  In
this case, we write $L_{T}^{(n)}$ to indicate its degree.

An interesting subalgebra of $\qs$ is the subspace $\Pi$ of \emph{peak functions}, which can be
defined as follows.  Let $\bc$ and $\bd$ be noncommuting indeterminates of degree 1 and 2,
respectively.  We let $w$ be an arbitrary word in the letters $\bc$ and $\bd$.  If
$$w=\bc^{n_{1}}\bd \bc^{n_{2}}\bd \cdots \bc^{n_{k}}\bd\bc^{n_{0}},$$
($n_{0}, \ldots , n_{k} \geq 0$)
then let $m_{j}= \deg(\bc^{n_{1}}\bd\bc^{n_{2}}\bd\cdots \bc^{n_{j}}\bd)$, $j=1,\dots,k$.
Define $\I^{w}$ to be the family consisting of the $k$ 2-element subsets $\{m_{j}-1,m_{j}\}$, $j=1,\dots,k$,
and
$$b[\I^{w}]\overset{\rm def}{=}
 \{ T \subseteq [n] : S\cap T \ne \emptyset, \hbox{\rm ~for all~} S\in \I^{w}\},$$
where n=$\deg(w)$.  Finally, define
\begin{equation}
\Theta_{w} \overset{\rm def}{=} \sum_{ T \in B[\I^{w}] } L_{T}^{(n+1)},
\label{thetaw}
\end{equation}
where $B[\I^{w}] \overset{\rm def}{=} \{ T \subseteq [n]: \; T,\overline{T} \in b[\I^{w}] \}$.
 Note that deg$(\Theta _{w})=\deg(w)+1$.
 Here, if $w={\bf 1}$, the empty $\cd$-word, then $\Theta_{\bf 1} = L_{\emptyset}^{(1)}=L_{1}$.
 We can define $\Pi$ to be the linear subspace of $\qs$ spanned
by $1$ and all the $\Theta_{w}$, as $w$ ranges over all $\cd$ words.
Again, $\Pi$ is a graded algebra with the inherited grading; we denote by 
$\Pi_{i}\overset{\rm def}{=} \Pi \cap \QQ_{i}$
its i$^{th}$ homogeneous part.  See \cite{BHV} and \cite{Ste} for details.  There the basis
element $\Theta_{w}$ in \eqref{thetaw} is replaced by
$\Theta_{w}^{(st)}=2^{\abs{w}_{\bd}+1}\Theta_{w}$, where
we denote the degree of $w$ by $\abs{w}$ and extend this notation to let $\abs{w}_{\bd}$
denote the number of $\bd$'s in $w$.

We summarize here the basics of the use of quasisymmetric functions in the theory
of flag $f$-vectors of graded posets and, in particular, Eulerian posets.  For a finite
graded poset $Q$, with rank function $\rk(\cdot)$, we define the formal power series
\begin{equation}
F(Q) \overset{\rm def}{=} \sum_{\hat{0} = u_0  \le \cdots \le u_{k-1}<u_{k} = \hat{1}}
x_1^{\rk(u_0,u_1)}  x_2^{\rk(u_1,u_2)} \cdots x_k^{\rk(u_{k-1},u_k)},
\label{Eqsf}
\end{equation}
where the sum is over all \emph{multichains} in $Q$
whose last two elements are different
and $\rk(x,y)=\rk(y)-\rk(x)$.
For general properties of $F(Q)$, see \cite{Eh} and \cite{BHV}.  In particular, we have
the following

\begin{proposition}For a graded poset $Q$,
\begin{enumerate}
\item[1.]
$F(Q) \in \qs$ and is homogeneous of degree $\rk(Q)$,
\item[2.]
$F(Q_1 \times Q_2) = F(Q_1) F(Q_2)$,
\item[3.]
$F(Q) = \sum_{\alpha} f_\alpha M_\alpha = \sum_{\alpha} h_\alpha L_\alpha$, where
$f_{\alpha}$ and $h_{\alpha}$ are the flag $f$ and flag $h$-vectors, respectively, of $Q$, and
\item[4.]
If $Q$ is Eulerian, then $F(Q) \in \Pi$; in fact $F(Q) = \sum_w [w]_{Q}\ \Theta_w$, where
 $[w]_{Q}$ denotes the coefficient of $w$ in the $\cd$-index of $Q$.
\end{enumerate}
\label{Fprops}
\end{proposition}

The last statement of Proposition \ref{Fprops}, which follows from  Proposition 2.2 and Theorem 2.1
of \cite{BHV}, can be used as the definition of the $\cd$-index for an Eulerian poset $Q$.  Formally,
the $\cd$-index is the homogeneous noncommutative polynomial,
$\psi_{Q} = \psi_{Q}(\bc,\bd) = \sum_{w}[w]_{Q} \thinspace w$ in $\bc$ and $\bd$, where the sum
is over all $\cd$-words of degree $\rk(Q)-1$; see \cite{bk}.

\section{The $R$-quasisymmetric function of a Bruhat interval and the complete $\cd$-index\label{sect2}}

Since a Bruhat interval $[u,v]$ is an Eulerian poset, it has a homogeneous $\cd$-index
$\psi_{u,v}  \overset{\rm def}{=} \psi_{[u,v]}$ and quasisymmetric function
$F(u,v)  \overset{\rm def}{=} F([u,v])$ defined as above.  
The polynomial $\psi_{u,v}$ has been studied explicitly by Reading \cite{Reading}.
In this section, we extend the definition of the $\cd$-index for Bruhat intervals to get a
nonhomogeneous polynomial, whose coefficients we later use to give a simple expression for
Kazhdan-Lusztig polynomials.

\subsection{The $R$-quasisymmetric function of a Bruhat interval}

We define a quasisymmetric function, analogous to the
power series (\ref{Eqsf}), making use of the polynomials $\widetilde{R}_{u,v}$ defined in
Proposition  \ref{Rtilde}.

Given $u,v \in W$, $u \leq v$, we define the \emph{R-quasisymmetric function}
$\widetilde{F}(u,v)$ by
\begin{equation}
\widetilde{F}(u,v) \overset{\rm def}{=} \sum_{u = u_0  \le \cdots \le u_{k-1} < u_k = v}
\widetilde{R}_{u_0u_1}(x_1)\widetilde{R}_{u_1u_2}(x_2)
\cdots \widetilde{R}_{u_{k-1}u_k}(x_k),
\label{Ftilde}
\end{equation}
where, again, the sum is over all multichains in $[u,v]$ whose last two elements are different.
Note that, by Proposition \ref{Rtilde} and the comments following Proposition \ref{2.2},
the leading term of each summand on the right hand side of (\ref{Ftilde}) is the
corresponding monomial on the right hand side of (\ref{Eqsf}).

An alternative description of $\widetilde{F}(u,v)$ in terms of \emph{chains} is as follows.
We omit the straightforward verification.

\begin{proposition}
Let $u,v \in W$, $u<v$. Then
\[ \widetilde{F}(u,v)=\sum_{k \geq 1} \; \sum_{u=u_{0}<u_{1}< \cdots <u_{k}=v} \;
\sum_{1 \leq i_{1}<i_{2}< \cdots <i_{k}} \; \prod_{j=1}^{k} \widetilde{R}_{u_{j-1},u_{j}}(x_{i_{j}}). \]
\label{chainform}
\end{proposition}

We will see that $\widetilde{F}(u,v)$ shares many of the properties of $F(Q)$ and will serve to
define an extension of the $\cd$-index for Bruhat intervals.  In particular, we have the following

\begin{theorem}
For any $u,v \in W$, $u \leq v$,
\begin{enumerate}
\item[1.]
$\widetilde{F}(u,v) = \sum_{\alpha} c_\alpha(u,v)\thinspace M_\alpha
 = \sum_{\alpha} b_\alpha(u,v)\thinspace L_\alpha$, and
\item[2.]
$\widetilde{F}(u,v) \in \Pi$, in fact
$$\widetilde{F}(u,v) \in \Pi_{l(u,v)} \oplus \Pi_{l(u,v)-2}\oplus \Pi_{l(u,v)-4}\oplus \cdots.$$
\end{enumerate}
\label{Ftildethm}
\end{theorem}

\begin{proof} 
To prove part 1,  we have, using Proposition \ref{2.LMS},
\begin{equation*}
\begin{split}
\widetilde{F}(u,v) &= \sum_{u = u_0 \le \cdots \le u_{k-1} < u_k = v}
\widetilde{R}_{u_0u_1}(x_1)\widetilde{R}_{u_1u_2}(x_2)
\cdots \widetilde{R}_{u_{k-1}u_k}(x_k) \cr
&= \sum_{u = u_0 \le \cdots \le u_{k-1} < u_k = v} \; \;
\prod _{j=1}^{k} \;  \sum_{\alpha _{j} \geq 0}\left( x_{j}^{\alpha _{j}}\;[q^{\alpha _{j}}]
( \widetilde{R}_{u_{j-1},u_{j}}) \right) \cr
&=  \sum_{u = u_0 \le \cdots \le u_{k-1} < u_k = v}\ \sum_{\alpha\in \N^{k}}\
\prod_{j=1}^{k}\left( x_{j}^{\alpha_{j}}\; [q^{\alpha_{j}}](\widetilde{R}_{u_{j-1},u_{j}})\right) \cr
&=  \sum_{k \geq 1}\ \sum_{\{ \alpha \in \N^{k}:  \alpha _{k}>0\}} x_{1}^{\alpha_{1}} \cdots
x_{k}^{\alpha _{k}} \sum_{u=u_{0} \le \cdots \le u_{k-1} <u_{k}=v}\ \prod_{j=1}^{k}[q^{\alpha_{j}}]
\left( \widetilde{R}_{u_{j-1},u_{j}} \right) \cr
&= \sum_{k \geq 1}\ \sum_{\{ \alpha \in \N^{k}:  \alpha _{k}>0\}} x_{1}^{\alpha_{1}} \cdots
x_{k}^{\alpha _{k}} \sum_{u=u_{0} < \cdots <u_{l(\alpha^{+})}=v}\ \prod_{j=1}^{l(\alpha^{+})}
[q^{\alpha_{j}^{+}}]\left( \widetilde{R}
_{u_{j-1},u_{j}} \right) \cr
&= \sum_{k \geq 1}\ \sum_{\{ \alpha \in \N^{k}:  \alpha _{k}>0\}} x_{1}^{\alpha_{1}} \cdots
x_{k}^{\alpha _{k}} \ c_{\alpha^{+}} (u,v) \cr
&= \sum_{\beta \in C} c_{\beta}(u,v)\ M_{\beta},
\cr
\end{split}
\end{equation*}
where  we have
used the facts that $[q^{0}](\widetilde{R}_{x,y})=\delta _{x,y}$ for all $x,y \in W$, $x \le y$,
and $\widetilde{R}_{x,x}=1$.
Here $\alpha^{+}$
is the composition obtained by taking only positive entries of $\alpha$ in order.
The second equality in part 1 follows from (\ref{2.IE.CB}) and (\ref{MtoL}).

Part 2 now follows easily  from
\cite[Theorem 8.4]{BrJAMS} and \cite[Proposition 1.3]{BHV}.
The last assertion follows since the $c_{\alpha}(u,v)$ count certain directed paths 
from $u$ to $v$ of length $|\alpha |$ in the Bruhat graph
$B(W,S)$, and all of these must have length $\equiv l(u,v)(\text{mod }  2)$.
\end{proof}

As a consequence of Theorem \ref{Ftildethm}, we can express any $\widetilde{F}(u,v)$ in
terms of the basis $\Theta_w $ for $\Pi$.
\begin{corollary}
\label{2.3}
For any  $u,v \in W$, $u \leq v$, we can write
\begin{equation*}
\widetilde{F}(u,v) = \sum_w [w]_{u,v}\ \Theta_w.
\end{equation*}
\label{extendcd}
\end{corollary}

 Note that, by Theorem \ref{Ftildethm}, the coefficients
$[w]_{u,v}$ can be nonzero only when $\deg(w) = l(u,v)-1, l (u,v)-3, \dots.$
We find it convenient to define, for any $u,v \in W$,
$$\widetilde{\psi}_{u,v} =\widetilde{\psi}_{u,v}(\bc,\bd)
 = \sum_w [w]_{u,v} \, w,$$
  a nonhomogeneous, noncommutative polynomial
in the variables $\bc$ and $\bd$ that  will be
called the \emph{complete $\cd$-index}
of the Bruhat interval $[u,v]$. 

\begin{example}
Let $W=S_{4}$, $u=1234$ and $v=4231$. Choose the reflection ordering $(1,2)<(1,3)<(1,4)<
(2,3)<(2,4)<(3,4)$.  Then there are 73 Bruhat paths from $u$ to $v$ and from (\ref{balpha})
one can compute that \medskip \\
$b_{5}(u,v)=b_{1,1,1,1,1}(u,v)=b_{1}(u,v)=1,\\
b_{3}(u,v)=b_{2,1}(u,v)=b_{1,2}(u,v)=b_{1,1,1}(u,v)=b_{1,4}(u,v)=b_{2,1,1,1}(u,v)=2,\\
b_{4,1}(u,v)=b_{3,1,1}(u,v)=b_{1,1,3}(u,v)=b_{1,1,1,2}(u,v)=3,\\
b_{2,3}(u,v)=b_{1,2,1,1}(u,v)=4,\\
b_{3,2}(u,v)=b_{1,1,2,1}(u,v)=5,\\
b_{1,3,1}(u,v)=b_{2,1,2}(u,v)=6, ~\text{and}\\
b_{2,2,1}(u,v)=b_{1,2,2}(u,v)=8,$\medskip\\
so by Theorem
\ref{Ftildethm}
\begin{align*}
\widetilde{F}(1234,4231) &= L_{5}+3L_{4,1}+5L_{3,2}+4L_{2,3}+2L_{1,4} \\
&  \quad+ 3L_{3,1,1}+6L_{1,3,1}+3L_{1,1,3}+8 L_{2,2,1}+6 L_{2,1,2}+8L_{1,2,2} \\
&  \quad+ 2 L_{2,1,1,1}+4L_{1,2,1,1}+5L_{1,1,2,1}+3L_{1,1,1,2}+L_{1,1,1,1,1} \\
&   \quad+ 2 L_{3} +2L_{2,1}+2L_{1,2}+2L_{1,1,1}+L_{1} \\
& =  \Theta _{\bc^{4}} +  \Theta _{\bd\bc^{2}}+2 \Theta _{\bc\bd\bc}+2\Theta _{\bc^{2}\bd}+
2 \Theta _{\bd^{2}} +2 \Theta _{\bc^{2}} + \Theta _{\bf 1},
\end{align*}
and so $[\bd\bc^{2}]_{1234,4231}=1, [\bc\bd\bc]_{1234,4231}=2$, {\it etc.}
Thus we have
\begin{equation*}
\widetilde{\psi}_{1234,4231} = \bc^{4}+\bd\bc^{2}+2\bc\bd\bc +2\bc^{2}\bd+2\bd^{2}+2\bc^{2}+{\bf 1}.
\end{equation*}
\label{s4examp}
\end{example}

There is at least one case where (\ref{Eqsf}) and (\ref{Ftilde}) define the same element of $\qs$.
For $u,v \in W$ and $i \in \N$ let $\widetilde{F}_{i}(u,v)$ be the homogeneous component
of $\widetilde{F}(u,v)$ of degree $i$ (so $\widetilde{F}_{i}(u,v)=0$ unless $i \equiv l (u,v)
\pmod{2}$).
\begin{proposition}
Let $u,v \in W$, $u \leq v$. Then
$\widetilde{F}_{l (u,v)}(u,v)=F(u,v)$ and so 
$[w](\widetilde{\psi}_{u,v}) = [w](\psi_{u,v})$ when $\deg(w)= l(u,v)-1$. In particular,
 $ \widetilde{F}(u,v) = F(u,v)$ and $\widetilde{\psi}_{u,v} =  \psi_{u,v}$ whenever the
Bruhat interval $[u,v]$ is a lattice.
\label{FFtilde}
\end{proposition}
\begin{proof}
The first assertion follows immediately from (\ref{Eqsf}), (\ref{Ftilde}) and the
comments following Proposition \ref{2.2}. The second one follows
from (\ref{Eqsf}), (\ref{Ftilde}) and \cite[Corollary 6.5]{BrInv}.
\end{proof}

\subsection{Algebraic properties of $\widetilde{F}$}

We investigate some of the algebraic properties of the map $[u,v] \mapsto \widetilde{F}(u,v)$.
In \cite[Proposition 4.4]{Eh}, the map on posets, $P \mapsto F(P)$, is shown to be
a morphism of Hopf algebras.  The same holds for the map defined by $\widetilde{F}$.
To see this, note that the coproduct on posets used in \cite{Eh} restricts to one on
Bruhat intervals:
\begin{equation}
\Delta([u,v]) = \sum_{w\in [u,v]} [u,w] \otimes [w,v].
\label{coproduct}
\end{equation}
On $QSym$, we take the usual coproduct defined by
\begin{equation}
\Delta(M_{\gamma}) = \sum_{\alpha \cdot \beta =  \gamma} M_{\alpha}\otimes M_{\beta},
\label{coproduct2}
\end{equation}
where $\alpha\cdot \beta$ denotes concatenation of compositions.

We begin with an analog of Proposition \ref{Fprops}.2, which shows the map to be multiplicative
with respect to direct product.

\begin{proposition}
\label{mult}
Let $W_{1},W_{2}$ be two Coxeter groups, $W_{1} \times W_{2}$ be their direct product,
and $u_{1},v_{1} \in W_{1}$, $u_{2},v_{2} \in W_{2}$. Then
\[ \widetilde{F}((u_{1},u_{2}),(v_{1},v_{2}))=\widetilde{F}(u_{1},v_{1}) \;
\widetilde{F}(u_{2},v_{2}). \]
\label{product}
\end{proposition}
\begin{proof}
It follows from Proposition \ref{productform} and Proposition \ref{Rtilde}
 that for $u_1 \le z_{1} \le z_{2} \le v_{1}$ and $u_2 \le w_{1} \le w_{2} \le  v_{2}$
\begin{equation}
\widetilde{R}_{(z_{1},w_{1}),(z_{2},w_{2})}= \widetilde{R}_{z_{1},z_{2}}\,
\widetilde{R}_{w_{1},w_{2}}.
\label{Rtildeproduct}
\end{equation}
The proof then follows by the same limiting  argument used to prove multiplicativity
of the map $F$
\cite[Proposition 4.4]{Eh}, with the change that we now use the maps
 $\widetilde{\kappa}_{i}(u,v)= \widetilde{R}_{u,v}(x_{i})$, which are multiplicative by 
 \eqref{Rtildeproduct}.  We omit the details.
\end{proof}

That $\widetilde{F}$ is a coalgebra map is proved next.

\begin{proposition}
\label{comult}
The map $\widetilde{F}$ is a map of coalgebras, that is, for each $u,v \in W$,
$u \le v$
\begin{equation*}
\left( (\widetilde{F} \otimes \widetilde{F} )\circ  \Delta \right)  ([u,v]) = \Delta \left( \widetilde{F}(u,v)\right).
\end{equation*}
\end{proposition}
\begin{proof}
By Theorem \ref{Ftildethm}.1 and \eqref{coproduct2}, we can write
\begin{align*}
\Delta \left( \widetilde{F}(u,v)\right)
&= \sum_{\gamma}c_{\gamma}(u,v) \sum_{\alpha \cdot \beta = \gamma} M_{\alpha}\otimes M_{\beta}\\
&= \sum_{\alpha}\sum_{\beta} c_{\alpha \cdot \beta}(u,v)\;  M_{\alpha}\otimes M_{\beta}.
\end{align*}
Now, as an extension of \cite[Proposition 5.5.4]{BB} (and with virtually the same proof),
we have $c_{\alpha \cdot \beta}(u,v) = \sum_{w\in [u,v]} c_{\alpha}(u,w)\thinspace c_{\beta}(w,v)$.  Thus
\begin{align*}
\Delta \left( \widetilde{F}(u,v)\right)
&=\sum_{w\in [u,v]}\sum_{\alpha}\sum_{\beta}
c_{\alpha}(u,w) \thinspace c_{\beta}(w,v) M_{\alpha}\otimes M_{\beta}\\
&= \sum_{w\in [u,v]}\widetilde{F}(u,w)\otimes \widetilde{F}(w,v),
\end{align*}
completing the proof.
\end{proof}

\begin{remark}
Formally, let ${\mathcal C}$ be the graded vector space, over a field ${\bf k}$, spanned by
$1\in {\bf k}$ and all isomorphism classes of
Bruhat intervals $[u,v]$, $u<v$, in all Coxeter groups, where the elements of ${\bf k}$
have degree 0 and $\deg([u,v])= l(v)-l(u)$, and
where $[u,v] \cong [w,z]$ if there exists a directed graph isomorphism $f : [u,v] \rightarrow
[w,z]$ such that $D(\Delta) = D(f(\Delta))$ for all Bruhat paths $\Delta$ in $[u,v]$.
${\mathcal C}$ has \emph{multiplication} defined
via Cartesian product $[u_{1},v_{1}] \times [u_{2},v_{2}] = [(u_1,u_2),(v_1,v_2)]$, and \emph{comultiplication}
defined by $\Delta([u,v])= \sum_{w\in [u,v]} [u,w]\otimes[w,v]$, where $[u,u]$ is defined to be $1\in {\bf k}$.
We define a \emph{counit} $\epsilon$ on ${\mathcal C}$ by $\epsilon([u,v])=\delta_{0,l(v)-l(u)}$
for $u \le v$.
By \cite[Lemma 2.1]{Eh}, this defines ${\mathcal C}$ as a \emph{graded Hopf algebra}.  Propositions
\ref{mult} and \ref{comult} show that the linear map $\widetilde{F}:{\mathcal C} \to QSym$
induced by $\widetilde{F}$ is a morphism of graded Hopf algebras.
(Compare \cite[Proposition 4.4]{Eh};
note that our $\widetilde{F}$ should not be confused with the similarly named function in \cite{Eh}.)
\end{remark}

\subsection{Properties of $\widetilde{\psi} _{u,v}$}

Let $\ba$ and $\bb$ be two noncommuting indeterminates and $\Z \langle \ba, \bb
\rangle $ be the ring of noncommutative polynomials in $\ba$ and $\bb$ with coefficients in
$\Z$. Given $n \in \PP$ and a subset $T \subseteq [n]$ let $m_{T}^{(n)}$  be the
noncommutative monomial of degree $n$ in $\ba$ and $\bb$ whose $i$-th letter (from the left) is $\ba$
if $i \not \in T$ and $\bb$ if $i \in T$, for $i=1,\ldots , n$.
Given a Bruhat path $\Gamma = (u_{0}, u_{1},
\cdots , u_{k-1},u_{k} )$  of length $k$ we define its {\em weight}
to be $w(\Gamma ) \overset{\rm def}{=} m_{D(\Gamma )}^{(k-1)}$.
\begin{proposition}
\label{Psitilde}
Let $u,v \in W$, $u<v$. Then
\[ \widetilde{\psi} _{u,v} (\ba + \bb,\ba \bb + \bb \ba) =\sum_{\Gamma } w(\Gamma ) \]
where $\Gamma $ runs over all the Bruhat paths from $u$ to $v$.
\end{proposition}
\begin{proof}
Let $w=\bc^{n_{1}} \, \bd \, \bc^{n_{2}} \, \bd \cdots \bc^{n_{k}}\, \bd \, \bc^{n_{0}}$ be a
$\cd$-word ($n_{0},\ldots , n_{k} \in \N$, $k \geq 0$). It then follows immediately
from the definition of $B[\I^{w}]$ that
\begin{equation}
\label{cd2ab}
(\ba + \bb)^{n_{1}} \, (\ba \bb + \bb \ba) \cdots (\ba + \bb)^{n_{k}} (\ba \bb + \bb \ba)(\ba + \bb)^{n_{0}}
= \sum_{T \in B[\I^{w}]} m_{T}^{(|w|)} .
\end{equation}
On the other hand, by Theorem 2.2 and Corollary 2.3,
\begin{equation}
\label{2.9.1}
\widetilde{F}(u,v) = \sum _{\alpha } b_{\alpha }(u,v) \, L_{\alpha } = \sum_{w}
[w]_{u,v} \, \Theta _{w} ,
\end{equation}
using (\ref{thetaw}) and equating coefficients of $L_{T}^{(n+1)}$ in (\ref{2.9.1})
we obtain that
\begin{equation}
\label{w2b}
b_{T}^{(n+1)}(u,v) = \sum_{\{ w: \; |w|=n, \; T \in B[I^{w}]\}} [w]_{u,v} ,
\end{equation}
for all $n \in \N$ and $T \subseteq [n]$, where  as in \eqref{balpha}
$$b_{T}^{(n+1)}(u,v) \overset{\rm def}{=}
|\{ \Gamma \in B_{n+1}(u,v) : \; D(\Gamma )=T \}|.$$
Therefore, by (\ref{cd2ab}) and (\ref{w2b}),
\begin{eqnarray*}
\widetilde{\psi}_{u,v}(\ba + \bb,\ba \bb + \bb \ba) & = & \sum_{w} [w]_{u,v} \, \sum_{T \in B_[\I^{w}]} m_{T}^{(|w|)} \\
& = & \sum_{n \geq 0} \, \sum_{T \subseteq [n]} \left( \sum_{\{ w: \; |w| =n, \, \, T \in B [I
^{w}] \} } [w]_{u,v} \right) m_{T}^{(n)} \\
& = & \sum _{n \geq 0} \, \sum_{T \subseteq [n] } b_{T}^{(n+1)} (u,v) \, m_{T}^{(n)} \\
& = & \sum _{\Gamma } w(\Gamma ) ,
\end{eqnarray*}
as desired.
 \end{proof}
 
 From this we see that we can obtain the polynomial $\widetilde{R}_{u,v}$ directly
 from $\widetilde{\psi} _{u,v}$.
 \begin{corollary}
 \label{psiRtilde}
For $u<v$, $\widetilde{R}_{u,v}(q) = q\; \widetilde{\psi}_{u,v}(q,0)$.
\end{corollary}
\begin{proof}
From \eqref{w2b}, one obtains $[\bc^{n}]_{u,v}= b^{(n+1)}_{\emptyset}(u,v)$ for
any $n$.  Thus
\begin{align*}
q\; \widetilde{\psi}_{u,v}(q,0) &= \sum_{n \geq 0}\; [\bc^{n}]_{u,v} \;q^{n+1}\\
		& = \sum_{\Gamma \in B(u,v): D(\Gamma)=\emptyset} q^{l(\Gamma)} = \widetilde{R}_{u,v}(q),
\end{align*}
the last equality being \cite[Theorem 5.3.4]{BjBr}.
\end{proof}
Thus, referring to Example \ref{s4examp}, we have $\widetilde{R}_{1234,4231}(q)=q^{5}+2q^{3}+q$.

As in \cite{ER}, we note that $\Z \langle \ba,\bb \rangle $ has a comultiplication defined
for monomials by
\begin{equation}
\label{coprod}
\Delta ' (a_{1} \cdots a_{n} ) \overset{\rm def}{=} \sum _{i=1}^{n} a_{1} \cdots a_{i-1}
\otimes a_{i+1} \cdots a_{n} ,
\end{equation}
 and extended linearly (where $\Delta ' (1) \overset{\rm def}{=} 0 \otimes 0$).  Similarly,
we can define a second coproduct on Bruhat intervals by
$$\Delta'([u,v])=\sum_{u<z<v}[u,z] \otimes [z,v].$$
(Note that this is \emph{not} the coproduct given in \eqref{coproduct}.)
We define, for convenience,
\begin{equation*}
\widetilde{\phi} _{u,v} (\ba,\bb) \overset{\rm def}{=} \widetilde{\psi}_{u,v}(\ba + \bb,\ba \bb + \bb \ba) ,
\end{equation*}
for all $u,v \in W$, $u<v$.  Then with the coalgebra structures just defined on Bruhat intervals
and $\Z \langle \ba,\bb \rangle $, we have that  $\widetilde{\phi}$ is a map of coalgebras.

\begin{proposition}
Let $u,v \in W$, $u<v$. Then
\[ \Delta ' (\widetilde{\phi }_{u,v}) = \sum _{u<z<v} \widetilde{\phi}_{u,z} \otimes \widetilde{\phi}
_{z,v} . \]
\end{proposition}
\begin{proof}
Given a Bruhat path $\Gamma = (u_{0}, u_{1}, \ldots
, u_{k} )$ we let $\Gamma _{(i)} \overset{\rm def}{=} (u_{0} , u_{1},
 \ldots , u_{i})$ and $\Gamma ^{(i)} \overset{\rm def}{=} (u_{i}
, u_{i+1} , \ldots , u_{k})$, for $i=1,\ldots , k-1$.
By Proposition \ref{Psitilde} and (\ref{coprod}) we have that
\begin{eqnarray*}
\Delta ' (\widetilde{\phi}_{u,v}) & = & \sum_{\Gamma } \Delta ' (w(\Gamma )) \\
& = & \sum _{\Gamma } \, \sum _{i=1}^{l(\Gamma )-1} w(\Gamma _{(i)}) \otimes w(\Gamma ^{(i)}) \\
& = &  \sum _{u<z <v} \, \sum_{\Gamma_{1} \in B(u,z) } \, \, \sum_{ \Gamma_{2} \in B (z,v)}
w(\Gamma_{1}) \otimes w(\Gamma_{2}) \\
& = & \sum _{u <z<v} \widetilde{\phi}_{u,z} \otimes \widetilde{\phi}_{z,v} ,
\end{eqnarray*}
as desired.
 \end{proof}

If we let $\bc=\ba + \bb$ and $\bd = \ba\bb + \bb\ba$, then as in \cite{ER}, the subalgebra
$\Z \langle \bc,\bd \rangle $ of $\Z \langle \ba,\bb \rangle $ is closed under $\Delta'$, since
$\Delta'(\bc)= 2(1\otimes 1)$ and $\Delta'(\bd) = \bc\otimes 1 + 1 \otimes \bc$ ($\Delta'$ acts
as a derivation on $\Z \langle \ba,\bb \rangle $).  Thus $\widetilde{\psi}$ is a map of coalgebras:
$$\Delta'(\widetilde{\psi}_{u,v})= \sum_{u<z<v} \widetilde{\psi}_{u,z} \otimes \widetilde{\psi}_{z,v}.$$

We note finally that the multiplicative structure used in \cite{ER}, the poset \emph{join}, does not
carry over to Bruhat intervals.  In most cases the join of two Bruhat intervals cannot be a Bruhat
interval by \cite[Theorem 3.2]{BrCaMa}.


\section{Kazhdan-Lusztig polynomials and the complete $\cd$-index\label{sect3}}

In this section we relate the Kazhdan-Lusztig
polynomial of any Bruhat interval in an arbitrary Coxeter group to its complete
$\cd$-index.

Recall the definition of the polynomials
$\Psi_\alpha$ and $\Upsilon_\alpha$ from \S 1.1.
Consider the map ${\mathcal K} : \qs \rightarrow \Z[q^{1/2},q^{-1/2}]$ defined by
$\K(M_\alpha) \overset{\rm def}{=}
(-1)^{l(\alpha)} q^{- \frac{|\alpha|}{2}} \Psi_\alpha$
or, equivalently,
$\K(L_\alpha) \overset{\rm def}{=} q^{- \frac{|\alpha|}{2}} \Upsilon_\alpha.$
Then by Theorem \ref{Ftildethm} we may rephrase Theorem \ref{1.8} in the following
way.
\begin{proposition}
\label{3.1}
Let $u,v \in W$, $u<v$.  Then
\begin{equation}
\K(\widetilde{F}(u,v)) =
q^{\frac{-l(u,v)}{2}} P_{u,v}(q) - 
q^{\frac{l(u,v)}{2}} P_{u,v} \left( \frac{1}{q} \right) . 
\end{equation}
\end{proposition}

If we define $\Xi_w \overset{\rm def}{=} \K(\Theta_w)$, then the following is
immediate from Corollary \ref{2.3} and Proposition \ref{3.1}.

\begin{corollary}
Let $u,v \in W$, $u<v$.  Then 
$$q^{\frac{-l(u,v)}{2}} P_{u,v}(q) - 
q^{\frac{l(u,v)}{2}} P_{u,v}(1/q) =
\sum_w [w]_{u,v}~\Xi_w.$$
\label{cdKL}
\end{corollary}

We next give an explicit description of $\Xi_{w}$ in terms of $w$.
Let $C_i \overset{\rm def}{=} \frac{1}{2i+1}\binom{2i+1} {i}$ for $i \in \N$
be the $i$-th Catalan number, and set $C_i = 0$ if $i \notin \N$.
Let
\begin{equation}
B_k(q) \overset{\rm def}{=} \sum_{i=0}^{\lfloor k/2 \rfloor} \frac{k+1-2i}{k+1}
\binom{k+1}{i}q^i.
\label{ballotpoly}
\end{equation}
We call $B_k(q)$ the $k$-th
\emph{ballot polynomial} since it is closely related to ballot problems
(see, {\it e.g.}, \cite[\S III.1, p. 73]{Fe}).

Let $a,b \in \Z$, $a \leq b$. By a {\em lattice path} on $[a,b]$
 we mean a function $\Gamma : [a,b] \rightarrow \Z$ such that
 $\Gamma (a) =0$ and
 \[ | \Gamma (i+1) - \Gamma (i)  |=1 \]
 for all $ i \in [a,b-1]$.  Given such a lattice path $\Gamma $ we let
 \[ N(\Gamma ) \overset{\rm def}{=} \{ i \in [a+1,b-1]: \Gamma (i)<0 \} ,
 \]
 \[ d_{+} (\Gamma ) \overset{\rm def}{=} | \{ i \in [a,b-1]: \Gamma
 (i+1) - \Gamma (i) =1 \} | , \]
  $l(\Gamma )  \overset{\rm def}{=} b-a$, and
  $ \Gamma _{\geq 0} \overset{\rm def}{=} l(\Gamma )-1-|N(\Gamma  )|$.
  We call $l(\Gamma )$ the {\em length} of
 $\Gamma $. Note that $b \not \in N(\Gamma )$ and that
 \begin{equation}
 \label{6.0.5}
 d_{+}(\Gamma )  = \frac{\Gamma (b)+b-a}{2} .
 \end{equation}
 We denote by ${\cal L}(n)$ the set of all lattice paths on $[0,n]$.
From \cite[Theorem 6.1]{BrJAMS} we have that for a composition $\alpha$,
\begin{equation}
\label{upsiloneq}
\Upsilon_{\alpha}(q) = (-1)^{|\alpha|-l(\alpha)}\sum_{\{\Gamma\in {\cal L}(|\alpha|): \, co(N(\Gamma))=\alpha^{*}\}} (-q)^{d_{+}(\Gamma)}.
\end{equation}
For a set $T$, we denote by $T^{*}$ the set for which
$co(T^{*})= co(T)^{*}$.

The next result, along with the previous one, shows how one can compute the Kazhdan-Lusztig
polynomial of any pair of elements $u,v \in W$ from the complete $\cd$-index of
the Bruhat interval $[u,v]$.

\begin{theorem}
Let $w=\bc^{n_{0}}\bd\bc^{n_{1}}\bd \cdots\bd\bc^{n_{k}}$,
$k \geq 0$, $n_{0} , \ldots ,n_{k} \geq 0$. Then
$$\Xi_{w} =   
(-1)^{k+\frac{\abs{w}-n_0}{2}}  \,
\left( q^{-\frac{n_0+1}{2}} B_{n_{0}}(-q)-q^{\frac{n_{0}+1}{2}} \, B_{n_{0}} \left( \frac{-1}{q} \right) \right)
 \, \prod_{j=1}^{k} C_{\frac{n_{j}}{2}}. $$
In particular, $\Xi_{w} = 0$ unless 
$n_{1} \equiv n_{2} \equiv \cdots \equiv n_{k} \equiv 0 \pmod{2}$.
\label{Xi}
\end{theorem}

\begin{proof}
Let $n = \abs{w}$, and for a $\cd$-word $w=w_{1}\cdots w_{n}$,
let $w^{\ast}\overset{\rm def}{=}w_{n}\cdots w_{1}$.
Note that $({\cal I}^w)^{\ast}={\cal I}^{w^{\ast}}$.
Then, by (\ref{thetaw}) and \eqref{upsiloneq} we have that
\begin{align}
q^{\frac{n+1}{2}} \Xi_w
 = &  \sum_{T \in B[{\cal I}^{w}]} \Upsilon_{co(T)} \nonumber  \\
 = & \sum_{T \in B[{\cal I}^{w}]}  \, \, \sum_{\{ \Gamma \in {\cal L}(n+1)
\, : \; N(\Gamma )^{\ast}=T \}} (-1)^{n+1-|N(\Gamma )|-1}(-q)^{d_{+}(\Gamma ) } \nonumber \\
\label{3.2}
 = & \sum_{\{ \Gamma \in {\cal L}(n+1)
\, , \; N(\Gamma ), \overline{N(\Gamma)}\in b({\cal I}^{w^{\ast}})\}} (-1)^{\Gamma _{\geq 0}}
\, (-q)^{d_{+}(\Gamma )} .
\end{align}
Now, ${\cal I}^{w^{\ast}}=\{ \{ p_{1}-1,p_1 \} , \{ p_2 -1, p_2\}, \ldots ,
\{ p_k -1, p_k \} \}$ where
$p_{j}\overset{\rm def}{=} n_{k+1-j}+\cdots +n_{k}+2j$ for $j \in [k]$, and clearly
$\overline{N(\Gamma )}= \{ i \in [n]: \; \Gamma (i)
\geq 0 \}$. Therefore $\Gamma \in {\cal L}(n+1)$ is such that $N(\Gamma ) \cap \{
p_{j}-1,p_{j} \} \neq \emptyset$
 and $ \overline{N(\Gamma )} \cap \{ p_{j}-1,p_{j} \} \neq \emptyset$ for all $j \in
[k]$ if and only if either $\Gamma (p_{j}-1)=\Gamma(p_{j})-1=-1$ or $\Gamma (p_{j}-1)
=\Gamma (p_{j})+1=0$, for all $j \in [k]$. Since
$\Gamma$ is a lattice path this happens if and only if
\begin{equation}
\label{3.3}
\Gamma (p_{j}-1) = \left\{ \begin{array}{ll}
\Gamma(p_{j})-1=-1, & \mbox{if $p_{j} \equiv 0 \pmod{2}$,} \\
\Gamma (p_{j})+1=0, & \mbox{if $p_{j} \equiv 1 \pmod{2}$,}
\end{array} \right.
\end{equation}
for all $j \in [k]$.

Let, for brevity, ${\cal L}_{w}$ be the set of all lattice paths $\Gamma \in
{\cal L}(n+1)$ that satisfy condition (\ref{3.3}) so, by (\ref{3.2}),
\begin{equation}
\label{3.4}
q^{\frac{n+1}{2}} \Xi_{w} = \sum_{\Gamma \in {\cal L}_{w}}(-1)^{\Gamma _{\geq 0}} \, (-q)^{d_{+}
(\Gamma )} .
\end{equation}
We claim that
\begin{equation}
\label{3.5}
q^{\frac{n+1}{2}} \Xi _{w} = \sum _{\Gamma \in {\cal L}^{\ast}_{w}} (-1)^{\Gamma _{\geq 0}} \,
(-q)^{d{+}(\Gamma )}
\end{equation}
where 
\begin{equation}
\label{3.6}
{\cal L}^{\ast}_{w} \overset{\rm def}{=} \{ \Gamma \in {\cal L}_{w} : \; \Gamma (i) \neq 0 \;
 \mbox{ if } \; i \in [n+1] \setminus \{ p_{1}, \ldots , p_{k} \} \}.
\end{equation}

In fact, let $\Gamma \in {\cal L}_{w} \setminus {\cal L}^{\ast}_{w}$ and let
$$i_{0} \overset{\rm def}{=} \min \{ i \in [n+1] \setminus \{ p_{1}, \ldots , p_{k} \}
: \; \Gamma (i)=0 \}.$$
Note that $i_{0} \equiv 0 \pmod{2}$. Then $p_{r}<i_{0}<p_{r+1}$
 for some $r \in [0,k]$ (where $p_{0} \overset{\rm def}{=} 0$ and
$p_{k+1} \overset{\rm def}{=}n+2$). Since $\Gamma \in {\cal L}_{w}$ this implies
that $\Gamma (p_{r})=0$ (else, by (\ref{3.3}), $\Gamma (p_{r}-1)=0$, contradicting the choice of
$i_{0}$). Define $\Gamma': [0,n+1] \rightarrow \Z$ by
\[ \Gamma ' (i) \overset{\rm def}{=} \left\{ \begin{array}{ll}
-\Gamma (i), & \mbox{if $p_{r} \leq i \leq i_{0}$,} \\
\Gamma (i), & \mbox{otherwise,}
\end{array} \right. \]
for $i \in [0,n+1]$. Then $\Gamma ' \in {\cal L}_{w} \setminus {\cal L}^{\ast}_{w}$,
$d_{+}(\Gamma ')=d_{+}(\Gamma )$ and
$(-1)^{\Gamma '_{\geq 0}}= -(-1)^{\Gamma_{\geq 0}}$. 
It is clear that
this map $\Gamma \mapsto \Gamma '$ is an involution of ${\cal L}_{w}
\setminus {\cal L}^{\ast}_{w}$ (since $i_{0}=$ min$\{ i \in [n+1] \setminus \{ p_{1},
\ldots , p_{k} \} : \; \Gamma ' (i)=0 \}$) and that it has no fixed points, so
$\sum _{\Gamma \in {\cal L}_{w} \setminus {\cal L}^{\ast}_{w}} (-1)^{\Gamma
_{\geq 0}} \, (-q)^{d_{+}(\Gamma)}=0$ and (\ref{3.5}) follows from (\ref{3.4}).

Now let $\Gamma \in {\cal L}^{\ast}_{w}$. Then, by (\ref{3.6}) and (\ref{3.3}),
\[ \Gamma (p_{j}-1)=\Gamma (p_{j})-1=-1 \]
and $p_{j} \equiv 0 \pmod{2}$ for all $j \in [k]$. This, again by (\ref{3.6}), implies that
$\Gamma (p_{j-1}+1)=-1$ for all $j \in [k]$ and that $\Gamma (i) \leq -1$ for all
$i \in [p_{k}] \setminus \{ p_{1}, \ldots , p_{k} \}$. Therefore, by (\ref{3.5}),
\begin{align*}
q^{\frac{n+1}{2}} \Xi _{w}  = &  \sum_{m=-n_{0}-1}^{n_{0}+1} \, \; \sum_{\{ \Gamma \in {\cal L}^{\ast}_{w} : \, \Gamma (n+1)=m \}} (-1)^{\Gamma _{\geq 0}} \, (-q)^{d_{+}(\Gamma)} \\
 = &
  \sum_{\{ m \in [n_{0}+1]: \;  m \equiv n_{0}+1 \pmod{2} \} } \;
 \prod_{j=1}^{k}C_{\frac{n_{j}}{2}} \,
\frac{m}{n_{0}+1} {\binom  {n_{0}+1} {\frac{n_{0}+1-m}{2}}} \\
 & \left( ( -1)^{k} \left( (-q)^{\frac{n+1+m}{2}} \, (-1)^{n_{0}}+(-q)^{\frac{n+1-m}{2}} \right) \right) \\
 = & \sum _{i=0}^{\lfloor \frac{n_{0}}{2} \rfloor}(-1)^{k} \prod_{j=1}^{k} C_{\frac{n_{j}}{2}}
\;
\frac{n_{0}+1-2i}{n_{0}+1} \binom  {n_{0}+1} {i} \\
 & \left( (-q)^{\frac{n+1+n_{0}+1-2i}{2}} \, (-1)^{n_{0}}+(-q)^{\frac{n+1-n_{0}-1+2i}{2}} \right) \\
 = & (-1)^{k} \, \prod_{j=1}^{k} C_{\frac{n_{j}}{2}}
 \, (-q)^{\frac{n-n_{0}}{2}}
\sum _{i=0}^{\lfloor \frac{n_{0}}{2} \rfloor } \frac{n_{0}+1-2i}{n_{0}+1}
\binom{n_{0}+1} {i} \\
 & \left( (-q)^{n_{0}+1-i} \, (-1)^{n_{0}}+(-q)^{i} \right)\\
 = & (-1)^{k} \, \prod_{j=1}^{k} C_{\frac{n_{j}}{2}} \, (-q)^{\frac{n-n_{0}}{2}}
\left( B_{n_{0}}(-q)-q^{n_{0}+1} \, B_{n_{0}}(-1/q) \right) ,
\end{align*}
where we have used (\ref{6.0.5}) and well known results on lattice path enumeration (see, {\it e.g.}, 
\cite[Ex.6.20]{ECII}).
\end{proof}

\section{Kazhdan-Lusztig polynomials and ballot polynomials\label{sect4}}

Using the results in the previous section, we derive the expansion,
in terms of the complete $\cd$-index, of the Kazhdan-Lusztig polynomials
with respect to a basis derived from the ballot polynomials \eqref{ballotpoly}.  This basis and its
relation to the Kazhdan-Lusztig polynomials was independently studied by Caselli \cite{Cas}, who
also considered its relation to the $R$-polynomials.

In what follows, we will consider a fixed Bruhat interval $[u,v]$; the dependence on the pair
$u,v$ will often be omitted.  We set $n = l(u,v)-1$.
For a $\cd$-word
$w=\bc^{n_{0}}\bd\bc^{n_{1}}\bd \cdots \bc^{n_{k-1}}\bd\bc^{n_{k}}$,
$k \geq 0$, $n_{0} , \ldots ,n_{k} \geq 0$,
define $C_{w}\overset{\rm def}{=}  \prod_{j=1}^{k} C_{\frac{n_{j}}{2}}$.   
If $w=\bc^{m}$, we take $C_{w}=1$.
Note that since $C_{i}=0$ when $i\notin \N$, $C_{w}=0$ unless $w$ is an \emph{even}
$\cd$-word, that is, unless $n_{1}, \dots , n_{k}$ are all even.

It will be helpful to rework Corollary \ref{cdKL} to obtain an
expression for the Kazhdan-Lusztig polynomial in terms of alternating shifted ballot polynomials
$q^i B_{n-2i}(-q)$.  Note that $q^{i}B_{n-2i}(-q)$ has degree $\lfloor n/2 \rfloor$ and
lowest degree term $q^{i}$ with coefficient 1.  Thus the set of polynomials $q^{i}B_{n-2i}(-q)$,
$0\le i \le \lfloor \frac{n}{2} \rfloor$, form a basis for the space of polynomials of degree $\leq
\lfloor n/2 \rfloor$.  Depending on the parity of $n$, Caselli denoted this basis by $O_{j}$ or $E_{j}$
(see \cite[Theorem 6.5]{Cas}). 

We begin by deriving the expression for $P_{u,v}$ in terms of this basis as a function of the
complete $\cd$-index.

\begin{theorem}
Let $u,v \in W$, $u \leq v$. Then
$$P_{u,v}(q)= \sum_{i=0}^{\lfloor n/2 \rfloor}a_i\  q^i \thinspace B_{n-2i}(-q),$$
where
\begin{equation*}
a_i = a_{i}(u,v) = [\bc^{n-2i}]_{u,v} + \sum_{\bd w \hbox{\Small ~even}}
 (-1)^{\frac{\abs{w}}{2} + \abs{w}_\bd}
C_{\bd w} [\bc^{n-2i}\bd w]_{u,v}.
\end{equation*}
\label{KLballot}
\end{theorem}

\begin{proof}
Since $\deg(P_{u,v})\le \lfloor n/2 \rfloor$, it follows from Corollary \ref{cdKL} that
\begin{equation}
P_{u,v} = \sum_{w} [w]_{u,v}\ D_{ \frac{n}{2}}\left( q^{\frac{l(u,v)}{2}} \Xi_{w} \right).
\label{Puv1}
\end{equation}
If $w=\bc^{n_{0}}\bd\bc^{n_{1}}\bd\bc^{n_{2}} \cdots \bd\bc^{n_{k}}$, with
$\abs{w} \le n=l(u,v)-1$, then by Theorem \ref{Xi},
\begin{equation*}
q^{\frac{l(u,v)}{2}} \Xi _{w}= (-1)^{k+\frac{\abs{w}-n_0}{2}} \, C_{w} \,
\left[ q^{\frac{n-n_0}{2}} B_{n_{0}}(-q)-q^{\frac{n+n_{0}}{2}+1} \, B_{n_{0}} \left( \frac{-1}{q} \right) \right].
\end{equation*}
Thus, if $w$ is even and $[w]_{u,v} \neq 0$, then $n_{0} \equiv \abs{w} \equiv n \pmod{2}$ so
$$D_{ \frac{n}{2}}\left( q^{\frac{l(u,v)}{2}} \Xi_{w} \right) =  (-1)^{k+\frac{\abs{w}-n_0}{2}} \, C_{w} \,
 q^{\frac{n-n_{0}}{2}} B_{n_{0}}(-q),$$
and (\ref{Puv1}) becomes
\begin{equation}P_{u,v} = \sum_{w \text{~even}}
 (-1)^{\abs{w}_{\bd}+\frac{\abs{w}-h(w)}{2}}~ C_{w}~  [w]_{u,v} ~q^{\frac{n-h(w)}{2}}B_{h(w)}(-q),\end{equation}where 
 $h(\bc^{n_{0}}\bd\bc^{n_{1}}\bd \cdots \bc^{n_{k-1}}\bd\bc^{n_{k}}) \overset{\rm def}{=} n_{0}$
 denotes the
 \emph{head} of $w$.
 
 Collecting the terms corresponding to $\cd$ words $w$ with $h(w)=n-2i$ gives the statement
 of the theorem.
\end{proof}

\begin{remark}
Alternatively, we can write $a_{i}$ as a sum over subsets of $[i]$, namely
\begin{align*}
\sum_{\overset{S \subseteq [i]}{S=\{ i_{1}, \ldots ,i_{k} \}_{<}}}  (-1)^{i_{k}-k}
\, [\bc^{n-2i}\bd \bc^{2(i_{1}-1)}\bd&\bc^{2(i_{2}-i_{1}-1)}\cdots\\ \dots&\bd \bc^{2(i_{k}-i_{k-1}-1)}
]_{u,v} 
  \prod_{j=1}^{k} C_{i_{j}-i_{j-1}-1}.
\end{align*}
\end{remark}

\begin{remark}
We list the first few coefficients $a_{i}$ as functions of the complete $\cd$-index (recall
$n= l(u,v)-1$, and that we omit the dependence on $u,v$):
\begin{equation}\begin{split}
a_{0} &= [\bc^{n}] \cr
a_{1} &= [\bc^{n-2}\bd] +[\bc^{n-2}] \cr
a_{2} &= [\bc^{n-4}\bd^{2}] - [\bc^{n-4}\bd\bc^{2}] +[\bc^{n-4}\bd] + [\bc^{n-4}] \cr
a_{3} &= [\bc^{n-6}\bd^{3}] - [\bc^{n-6}\bd^{2}\bc^{2}] - [\bc^{n-6}\bd\bc^{2}\bd] +2[\bc^{n-6}\bd\bc^{4}] \cr
   &\phantom{xxxxxxxxxxx}+  [\bc^{n-6}\bd^{2}] - [\bc^{n-6}\bd\bc^{2}] +[\bc^{n-6}\bd] + [\bc^{n-6}].
\end{split}
\label{aicd}\end{equation}
\end{remark}

\begin{remark}
By Proposition \ref{FFtilde}, the degree $n$ terms of the complete \cd-index of the Bruhat
interval $[u,v]$ correspond to the ordinary \cd-index of the Eulerian poset $[u,v]$.
Restricting the formulas for the $a_{i}$ in Theorem \ref{KLballot} to
only the degree $n$ \cd-coefficients  yields the expression given by Bayer and Ehrenborg
\cite[Theorem 4.2]{BaEh} for the $g$-polynomial $g([u,v]^{*},q)$ in terms of the \cd-index
of the dual interval $[u,v]^{*}$.
This can be checked by comparing formulas -- for example the polynomial $Q_{k+1}(x)$ of
\cite{BaEh} is $B_{k}(-x)$ -- and recalling that
$[w]_{u,v} = [w^{*}]_{[u,v]^{*}}$.  One consequence
is that the difference $P_{u,v}(q)- g([u,v]^{*},q)$ is a function of the lower degree 
\cd-coefficients only.
\end{remark}

\begin{example}
Continuing with Example \ref{s4examp}, we have
\begin{equation*}
\widetilde{F}(1234,4231) 
 =  \Theta _{\bc^{4}} +  \Theta _{\bd\bc^{2}}+2 \Theta _{\bc\bd\bc}+2\Theta _{\bc^{2}\bd}+
2 \Theta _{\bd^{2}} +2 \Theta _{\bc^{2}} + \Theta _{\bf 1}
\end{equation*}
and so
\begin{align*}
P_{u,v}(q) &= a_{0}q^{0}B_{4}(-q)+a_{1}q^{1}B_{2}(-q)+a_{2}q^{2}B_{0}(-q)\\
&=[{\bc^{4}}] \left(1-3q+2q^{2}\right) + \left([{\bc^{2}\bd}]+ [{\bc^{2}}]\right) q(1-q)\\
&\phantom{a_{0}q^{0}B_{4}(-q)+a_{1}q^{1}B_{2}}
+\left([{\bd^{2}}]-[{\bd\bc^{2}}]+[{\bd}]+[{\bf 1}]\right)q^{2}\\
&=(1-q+q^{2})+(2q-q^{2})=1+q
\end{align*}
with $g([u,v]^{*},q)= 1-q+q^{2}$.
\end{example}

Next we give an example that shows that it is not possible to express $P_{u,v}$ as a function
of the ordinary (homogeneous) $\cd$-index $\psi_{u,v}$ alone.  That is, the Kazhdan-Lusztig
polynomial does not depend only on the flag $f$-vector of the Eulerian poset $[u,v]$.
\begin{example}
Let $W=S_{5}$ and consider the rank 6 Bruhat intervals $[12435,53142]$ and $[31254,53421]$.
One can compute
\begin{align*}
\widetilde{\psi}_{12435,53142} =  \bc^{5} &+ 6\bc\bd\bc^{2} + 6\bc^{2}\bd\bc + 3\bd\bc^{3} + 3\bc^{3}\bd 
	+  7\bc\bd^{2} + 7\bd^{2}\bc + 6\bd\bc\bd\\
	&   + \bc^{3} + 2\bd\bc + 2 \bc\bd    \\
\widetilde{\psi}_{31254,53421} =  \bc^{5} &+ 6\bc\bd\bc^{2} + 6\bc^{2}\bd\bc + 3\bd\bc^{3} + 3\bc^{3}\bd
	+  7\bc\bd^{2} + 7\bd^{2}\bc + 6\bd\bc\bd\\
	&   + 2\bc^{3} + 4\bd\bc + 4 \bc\bd
\end{align*}
while $P_{12435,53142}=1$ and $P_{31254,53421}= 1+q$.  Thus neither $\widetilde{\psi}_{u,v}$
nor $P_{u,v}$ is a function of $\psi_{u,v}$ alone.
\end{example}
It appears that no such example was previously known.  It can be checked that there is no such
pair of intervals in $S_{4}$.  There are many other examples in $S_{5}$, although none of them
involve lower intervals ({\it i.e.}, those with $u=e$).

\bigskip

Considering the dependence of the coefficients of the Kazhdan-Lusztig polynomials on
these coefficients $a_{i}$, let $P_{u,v}=p_0 + p_1 q + \cdots $.  The following is a direct
consequence of the definition of the $a_{i}$ as the coefficients of $P_{u,v}$ in the basis
$q^{i}\thinspace B_{n-2i}(-q)$.  Propositions \ref{pfroma} and \ref{afromp} as well as
Corollary \ref{nonneg} are all implicit in \cite[\S6]{Cas}, so proofs will be omitted here.

\begin{proposition}
For $j=0,\dots,\lfloor n/2 \rfloor$,
$$p_j = \sum_{i=0}^j (-1)^{j-i} \frac{n+1-2j}{n+1-2i}\binom{n+1-2i}{j-i}a_i. $$
\label{pfroma}
\end{proposition}

Theorem \ref{KLballot} and Proposition \ref{pfroma} allow us to derive the coefficients of
the Kazhdan-Lusztig polynomial for any Bruhat interval as a function of its complete $\cd$-index.
For example, we can read from Proposition \ref{pfroma} and (\ref{aicd}) that
\begin{equation}
\label{p1coeff}
\begin{split}
p_{1}=[q](P_{u,v})  = & a_{1}-(n-1)a_{0} \\
 = & [\bc^{n-2}\bd]+[\bc^{n-2}]-(n-1)[\bc^{n}].
 \end{split}
\end{equation}

The relations given in Proposition \ref{pfroma} are unitriangular so invertible.  The inverse
relations have a particularly simple nonnegative form.

\begin{proposition}
\label{afromp}
For $j=0,\dots,\lfloor n/2 \rfloor$,
$$a_j = \sum_{i=0}^j \binom{n-j-i}{n-2j}p_i.$$
\end{proposition}

\begin{corollary}
\label{nonneg}
Nonnegativity of $a_{i}(u,v)$, $i=1,\dots,k$ is implied by the nonnegativity of
$[q^{i}](P_{u,v})$, $i=1,\dots,k$.
\end{corollary}

\begin{remark}
In fact, \emph{any} nonnegative polynomial of degree $\leq \lfloor
\frac{n}{2} \rfloor$ has a nonnegative representation
in terms of the basis of alternating shifted ballot polynomials $q^i B_{n-2i}(-q)$.
\end{remark}

Given the conjectured nonnegativity of the coefficients of the Kazhdan-Lusztig polynomial
$P_{u,v}$ for any Bruhat interval $[u,v]$ \cite{K-L}, Caselli made
 the following conjecture.  Because of Theorem \ref{KLballot} (see also \eqref{aicd}), we can
 interpret it as a conjectured set of linear inequalities that must be satisfied by the complete
 $\cd$-index.

 \begin{conjecture}\cite[Conjecture 6.6]{Cas}
 For each Bruhat interval $[u,v]$ and for each $i=0,1,\dots , \lfloor n/2 \rfloor$,
 $n=l(u,v)-1$, we have $a_{i}(u,v)\ge 0$. 
  \end{conjecture}

\medskip
We conclude by noting the following consequence of Proposition \ref{afromp}.
\begin{corollary}
Let $u,v \in W$, $u<v$. Then 
$$a_{\lfloor \frac{l(u,v)-1}{2} \rfloor} =
\begin{cases}
P_{u,v}(1), &  l(u,v) \hbox{\rm ~odd} \cr
\frac{d}{dq}\left(q^{\frac{l(u,v)}{2}} P_{u,v}(1/q)\right)\Big|_{q=1},
& l(u,v) \hbox{\rm ~even.}
\end{cases}$$
\end{corollary}

\begin{remark}
An interesting question is whether, for all $w$, the quantities $[w]_{u,v}$ are
\emph{combinatorially invariant}, that is, they
depend only on the \emph{poset} structure of the interval $[u,v]$.  This is true for $\cd$-words
$w$ of degree $l(u,v)-1$ by Proposition \ref{FFtilde}.
By Theorem \ref{KLballot},  combinatorial invariance of $[w]_{u,v}$ for all $w$ implies
the conjectured combinatorial invariance of the Kazhdan-Lusztig polynomial $P_{u,v}$.
But the converse is also true: if $P_{u,v}$ is combinatorially invariant for all $u<v$, then
so is $[w]_{u,v}$ for all $w$ and all $u<v$.  This follows from Proposition \ref{2.LMS} and 
Theorem \ref{KLdefn} and the fact that $[w]_{u,v}$ depends on the $c_{\alpha}$, by
\eqref{2.IE.CB} and Proposition \ref{k2cd}.  Combinatorial invariance of $P_{u,v}$ is known
to hold in the case $u=e$  (see \cite{Mar}, \cite{BrCaMa} and \cite{Del}), and it follows
that the same is true for all $[w]_{e,v}$.
\end{remark}


\section{Other representations of the complete $\cd$-index \label{sect5}}

The formulas obtained in the previous sections for the Kazhdan-Lusztig
polynomials assume knowledge of the complete $\cd$-index.
In this section we show how one can explicitly compute the complete $\cd$-index
of a Bruhat interval $[u,v]$ in terms of the coefficients $b_{\alpha}(u,v)$ defined in
\S 1.1.

Let $u,v \in W$, $u<v$, and $n \in [l(u,v)-1]$, $ n \equiv l(u,v)-1 \pmod{2}$.
In this section it will be convenient to index the $b_{\alpha}$
by subsets instead of the corresponding compositions.
Unless otherwise explicitly stated, all subsets are in $[n]$.
We will omit to write the dependence on $u,v$ throughout the section.

We start with a general expression for $[w]$ in terms of \emph{sparse} $b_{S}$,
that is, those
$b_{S}, S\subseteq [n]$, where $S$ has no two successive elements and
$n \notin S$.  In terms of $b_{\alpha}$, this means that $\alpha_{i} > 1$
if $i>1$.
We use an expression for $[w]$ in terms of the \emph{sparse $k$-vector},
which is a reinversion
of the sparse $b_{S}$ defined by $k_{S}= \sum_{T\subseteq S} (-1)^{|S|-|T|}b_{T}$.
The following result is proved in exactly the same way as \cite[Proposition 7.1]{BE}.

\begin{proposition}
\label{k2cd}
Let $w = \bc^{n_1}\bd\bc^{n_2}\bd \cdots\bc^{n_k}\bd\bc^{n_{0}}$
 and define $m_0,\ldots,m_k$ by $m_0=1$ and
$m_i=m_{i-1}+n_i+2$.  Then
\begin{align} \label{kformula}
 [w] = \sum_{i_1,\ldots,i_k} (-1)^{(m_1-i_1) +(m_2-i_2) + \cdots
      + (m_k-i_k)}~k_{ \{ i_1 i_2 \cdots i_k \}},
  \end{align}
where the sum is over all $k$-tuples $(i_1,i_2,\ldots,i_k)$ such that
$m_{j-1} \leq i_j \leq m_j-2$.
\label{kprop}
\end{proposition}
The inversion of this relation, expressing $k_{S}$ as a sum of distinct $[w]$, is
the same as \cite[Proposition 2.3]{BHV}.  Using Proposition \ref{k2cd} one can obtain
from \eqref{p1coeff} the coefficient of the linear term of the Kazhdan-Lusztig
polynomial $P_{u,v}$ (compare \cite[Corollary 5.9]{BrInv} and \cite[Theorem D]{BjEk}).
\begin{corollary}  For any Bruhat interval $[u,v]$, the coefficient of the linear
term in $P_{u,v}$ is
$$
p_1=c_{n,1}(u,v) + c_{n-1}(u,v)-(n+1).
$$
\end{corollary}

We will show that if $[w]$ is expressed in terms of the sparse $b_S$,
the coefficients are still $\pm 1$, leading
to the possibility of a direct enumerative interpretation of $[w]$.
Using the notation of Proposition \ref{kprop}, let
$$A_{j} = \{\; i\;  :\;  m_{j-1}\le i\le m_{j}-2\}$$
be the range of $i_{j}$ above,
and let 
$$\T_{w}\overset{\rm def}{=} \{ T\subset \cup_{i\le k} A_{i} : |T\cap A_{j}| \le1~\forall j \}.$$
Then $k_{S}$ appears in \eqref{kformula} if and only if $S\in \T_{w}$ and $|S|=k$ ,
and by Proposition \ref{kprop} and the definition of $k_S$, we see that $[w]=\sum_{S\in\T_{w}}d_{S} \; b_{S}$ for
integers $d_{S}$.  We now prove that these integers can only be $\pm 1$ or $0$.
\begin{proposition}
\label{5.2}
If $w=\bc^{n_{1}}\bd\cdots \bc^{n_{k}}\bd\bc^{n_{0}}$,
then $[w]=\sum_{S\in \T_{w}} d_{S}\thinspace b_{S}$, where
$$d_{S} =
\begin{cases}
(-1)^{k-|S|}\prod_{j}\prod_{i\in S\cap A_{j}} (-1)^{m_{j}-i}, & \hbox{\rm ~if~} \abs{S\cap A_{i}}=1
\hbox{\rm ~for~} |A_{i}| \hbox{\rm ~ even,} \cr
0, & \hbox{\rm ~otherwise}.
\end{cases}$$
\end{proposition}

\begin{proof}
By \eqref{kformula} and the definition of $k_{S}$ we have
\begin{equation*}\begin{split}
d_{S}  &= \sum_{S\subseteq \{i_{1},\dots,i_{k}\} \in \T_{w}}(-1)^{k-|S|}
(-1)^{(m_1-i_1) +(m_2-i_2) + \cdots  + (m_k-i_k)}
\cr
&= (-1)^{k-|S|}\left[ \prod_{j: S\cap A_{j}\ne \emptyset}\ \prod_{i\in S\cap A_{j}} (-1)^{m_{j}-i}\right]\
\left[ \prod_{j: S\cap A_{j} = \emptyset}\ \left(\sum_{i\in A_{j}} (-1)^{m_{j}-i}\right)\right].
\cr
\end{split}\end{equation*}
But the sums in the last expression are 0 or 1 depending on whether $|A_{j}|$ are even or odd,
so the result follows.
\end{proof}

For example, for $w=\bc \bd^{2}$, $k=2$, $n_{1}=1$, $n_{2}=0$, $n_{0}=0$, $A_{1}=\{1,2\}$, $A_{2}=\{4\}$ and so
$$ [\bc \bd^{2}] =  b_{\{2,4\}} - b_{\{1,4\}} - b_{\{2\}} + b_{\{1\}}. $$
  Note that the same expression computes
 the coefficient $[w]=[\bc \bd^{2}\bc^{m}]$ for \emph{any} $m\ge 0$, where the subscripts are
 understood to be subsets of $\{1,\ldots,\deg w\}$.

Finally, we give simple proofs of two elementary identities between $\cd$-coefficients
and the numbers $b_{\alpha}$ using basic quasisymmetric identities.
The second of these is essentially in \cite{stancd}.

\begin{proposition}
\label{sums}
Let $u,v \in W$, $u<v$, and $n \in [l(u,v)-1]$, $ n \equiv l(u,v)-1 \pmod{2}$. Then:
\begin{enumerate}
\item
$\sum_{w} 2^{n-\abs{w}_{d}} [w]_{u,v} = c_{[n]}(u,v)$, the number of Bruhat paths from $u$ to $v$
of length $n+1$,
\item
$\sum_{w}[w]_{u,v} = b_{\{1,3,5,\dots \}}(u,v)$,
\end{enumerate}
where the sums are over all $\cd$-words of degree $n$.
\end{proposition}

\begin{proof}
By Theorem \ref{Ftildethm}
$$\widetilde{F}_{n+1}(u,v)=\sum_{S}c_{S}(u,v)\ M_{S}^{(n+1)}=\sum_{S}b_{S}(u,v)\ L_{S}^{(n+1)}=
\sum_{w} [w]_{u,v}\ \Theta_{w},$$
where the sums are over all subsets of $[n]$ and over all $\cd$-words of degree $n$
(recall that we denote by $\widetilde{F}_i$ the homogeneous component of
$\widetilde{F}$ of degree $i$).
We now use the relations
$$\Theta_{w}= \sum_{S\in b[\I_{w}]} \ 2^{|S|-|w|_{\bd}}\ M_{S}^{(n+1)} =
\sum_{T \in B[\I^w]} L_T^{(n+1)}$$
from \cite[(2.14)]{BHV} (divided, as before, by $2^{|w|_{\bd}+1}$) and \eqref{thetaw}.
Hence the coefficient of $M_{[n]}^{(n+1)}$ in $\Theta_{w}$ is  $2^{n-|w|_{\bd}}$, so we can conclude (1).
Similarly, the coefficient of $L_{\{ 1,3,5,\ldots \}}$ in $\Theta _{w}$ is $1$, so (2) follows.
\end{proof}


\section{A sign conjecture for the complete $\cd$-index\label{sect7}}
We conclude with an intriguing sign conjecture for the complete $\cd$-index and
prove one of its simple consequences in certain cases.
Computation on Bruhat intervals up to rank 7 occurring in symmetric groups
suggests the following.
\begin{conjecture}
\label{main}
Let $(W,S)$ be a Coxeter system and $u,v \in W$, $u<v$. Then
$$[w]_{u,v} \ge 0$$
for all $\cd$-words  $w$.

Since every Bruhat interval is shellable ({\it e.g.}, \cite[Theorem 2.7.5]{BjBr}) and Eulerian,
hence Gorenstein$^{*}$,
it follows from the recent work of Karu \cite{KaruComp} (see also \cite{EK})
that $[w]_{u,v} \ge 0$
whenever $\deg w = l(u,v)-1$ (the top degree for which $[w]_{u,v} \ne 0$).
\end{conjecture}

If Conjecture \ref{main} were true, then Proposition \ref{sums} would imply that
$$ 2^{n} [\bc^n]_{u,v} \leq c_{[n]}(u,v)$$
for all $u,v \in W$ and $n \geq 0$.
This is indeed true for finite Coxeter and affine Weyl groups, as we now show.
From now on, we assume that W is a finite Coxeter or affine Weyl group.

Given two Bruhat paths $\Delta  =(u, x ,v)$, $\Gamma =(u ,y , v) \in B_{2}(u,v)$
write $\Delta \leq \Gamma$ if $(xu^{-1},vx^{-1})$ is lexicographically smaller than $(yu^{-1},vy^{-1})$
(where $T$ is totally ordered by the chosen reflection ordering).
Let $m \overset{\rm def}{=} [t^{2}](\widetilde{R}_{u,v})$.
It then follows easily from  \cite[Theorem 5.3.4]{BjBr} that
 $|B_{2}(u,v)|=2m$ and that 
 $$|\{ \Delta \in B_{2}(u,v): \;
D(\Delta ) =\emptyset \} |=m=|\{ \Delta \in B_{2}(u,v) : \; D(\Delta ) = \{ 1 \} \} |.$$
For $\Gamma \in B_{2}(u,v)$ let
\[ r \overset{\rm def}{=} | \{ \Delta \in B_{2}(u,v): \; D(\Delta )=D(\Gamma ) , \; \Delta
\leq \Gamma \} |. \]
The {\em flip} of $\Gamma$ is the $r$-th Bruhat path (in the lexicographic ordering) in
 $\{ \Delta \in B_{2} (u,v): D(\Delta )
\neq D(\Gamma ) \} $. We denote this path by $\textit{flip}(\Gamma)$. This notion is a special case
of that of lexicographic correspondence first introduced in \cite{BreInc}. The following result is a special case of \cite[Corollary 2.3]{BreInc}.
\begin{proposition}
\label{flip}
Let $W$ be a finite Coxeter or affine Weyl group.  Further,
let $u,v \in W$, $u<v$, $(u,y,v) \in B_{2}(u,v)$ be such that $D((u,y,v))=
\emptyset$ and $(u,x,v) \overset{\rm def}{=} \textit{flip}((u,y,v))$. Then
$yu^{-1}<_{T}xu^{-1}$ and $vx^{-1}<_{T}vy^{-1}$.
\end{proposition}
Note that $\textit{flip}(\textit{flip}(\Gamma ))=\Gamma$. Given $\Delta =(u_{0},u_{1},
\ldots , u_{k} ) \in B_{k} (u, v)$ we let, for $i \in [k-1]$, $\textit{flip}_{i}
(\Delta )$ be the Bruhat path in $B_{k}(u,v)$ obtained by flipping $\Delta$ at its $i$-th element.
Namely, $\textit{flip}_{i}(\Delta ) \overset{\rm def}{=} (u_{0}, u_{1},
\ldots ,u_{i-1},$ $x , u_{i+1} ,
\ldots , u_{k-1},u_{k})$ where $(u_{i-1}, x,
u_{i+1})
\overset{\rm def}{=} \textit{flip} ((u_{i-1},u_{i},u_{i+1} ))$. Note that if $\{ i-1,i,i+1 \} \cap
D(\Delta ) = \emptyset $ then, by Proposition \ref{flip}, $\{ i-1,i,i+1 \} \cap D(\textit{flip} _{i}(\Delta ))=
\{ i \}$.
\begin{proposition}
\label{chaininject}
Let $W$ be a finite Coxeter or affine Weyl group and $u,v \in W$, $u<v$. Then
\[ [t^{k}] (\widetilde{R}_{u,v}) \leq \left\lfloor \frac{|B_{k}(u,v)|}{2^{k-1}} \right\rfloor \]
for all $k \geq 1$,
and so $ 2^{n} [\bc^n]_{u,v} \leq c_{[n]}(u,v)$ for $n \geq 0$.
\end{proposition}
\begin{proof}
Since $[t^{k}](\widetilde{R}_{u,v} ) \in \N$, we prove the equivalent statement that
\[ [t^{k}] (\widetilde{R}_{u,v}) 2^{k-1} \leq |B_{k} (u,v)| . \]
We do this by constructing an explicit injection.

Let $\Gamma =(u, u_{1} , \cdots
,u_{k-1} ,v ) \in B_{k}(u,v)$ be such that $D(\Gamma )=\emptyset$ and $S =
\{ s_{1}, \ldots , s_{r} \} _{<} \subseteq [k-1]$. We define
\[ \varphi (\Gamma ,S) \overset{\rm def}{=} \textit{flip} _{s_{r}}(\cdots (\textit{flip}_{s_{1}}(\Gamma )) \cdots
). \]
Clearly, $\varphi (\Gamma ,S ) \in B_{k}(u,v)$. Furthermore, note that by Proposition
\ref{flip}
\begin{eqnarray*}
s_{i} & =& \mbox{max} \; \{ D(\textit{flip}_{s_{i}}(\cdots (\textit{flip}_{s_{1}}(\Gamma )) \cdots
)) \} \\
& = & \mbox{max} \; \{ D (\textit{flip}_{s_{i+1}}(\cdots (\textit{flip}_{s_{r}}(\varphi (\Gamma ,S))) \cdots
)) \}
\end{eqnarray*}
for $i=1, \ldots , r$.Therefore $s_{r}, \ldots , s_{1}$ and hence $\Gamma$ are uniquely
recoverable from $\varphi (\Gamma , S)$, so $\varphi $ is an injection. The result now
follows from the  fact that $[t^{n+1}](\widetilde{R}_{u,v})=[\bc^{n}]_{u,v}$ by Corollary \ref{psiRtilde}.
\end{proof}

\begin{remark}
We note that Proposition \ref{chaininject} holds whenever Proposition \ref{flip} does,
and the latter has been conjectured to hold for all Coxeter groups (see \cite[p.\ 117]{Cel} and
\cite[p.\ 745]{BreInc}),
giving further evidence for Conjecture \ref{main}.  In fact, a proof of Proposition \ref{flip}
that holds for all Coxeter groups has been recently announced \cite{SAB}.
\end{remark}

\providecommand{\bysame}{\leavevmode\hbox
to3em{\hrulefill}\thinspace}

\end{document}